%% file: cxcrossed.tex
\documentclass[12pt]{amsart}
\usepackage{pb-diagram}
\usepackage{cmex12-2e}
\usepackage[6in]{amspapermacs}

\input cxcmacs

\begin{document}

\title[$C_0(X)$-Actions]{Crossed products 
by \boldmath $C_0(X)$-actions}

\author[Echterhoff]{Siegfried Echterhoff}
\address{Universit\"at-Gesamthochschule Paderborn \\
Fachbereich Mathematik-Informatik \\
Warburber Stra\ss e 100 \\
D-33095 Paderborn \\
Germany}
\email{echter@uni-paderborn.de}

\author[Williams]{Dana P. Williams}
\address{Department of Mathematics \\
Dartmouth College \\
Hanover, NH 03755-3551 \\
USA}
\email{dana.williams@dartmouth.edu}

\dedicatory{Dedicated to Professor E.~Kaniuth on the occasion of his
$\text{60}^{\text{th}}$ birthday} 

\date{16 June 1997}

\ifdraft
\include{cxcab}\else
\input cxcab
\fi

\maketitle

\ifdraft

\include{intro1}\else
\input intro1
\fi
\ifdraft
\include{bund-op1}\else
\input bund-op1
\fi
\ifdraft
\include{bund-op2}\else
\input bund-op2
\fi
\ifdraft
\include{inner}\else
\input inner
\fi
\ifdraft

\include{loc-fibre}\else
\input loc-fibre
\fi
\ifdraft
\include{loc-inner}\else
\input loc-inner
\fi
\ifdraft
\include{action}\else
\input action
\fi       
\includerefs

\end{document}

%% file: cxcmacs.tex
\makeatletter
\@addtoreset{footnote}{section}
\@addtoreset{footnote}{subsection}
\makeatother

\theoremstyle{definition}
\newtheorem{remark}[thm]{Remark}
\newtheorem{ex}[thm]{Example}

\theoremstyle{definition}
\newtheorem*{ga}{Standing Assumptions}
\newtheorem*{problem}{Open Question}
\theoremstyle{plain}
\newtheorem*{thmnn}{Theorem}

\DeclareMathSymbol{\rtimes}{\mathbin}{AMSb}{"6F}

\def\e{\xi}
\def\sec{\Gamma_0}
\def\s(#1){\sec(\e,#1)}
\def\hm{homomorphism}
\def\UM{UM}

\def\subab{\text{\normalfont ab}}
\def\gab{G_{\subab}}

\def\hgab{\widehat{G}_{\subab}}
\def\hlab{\widehat{L}_{\subab}}

\def\sheaf#1{\mathcal {#1}}

\def\shgab{\widehat{\sheaf G}_{\subab}}

\def\H{\mathcal{H}}

\def\<{\langle}
\def\>{\rangle}
\let\ipscriptstyle=\scriptscriptstyle
\def\lipsqueeze{{\mskip -3.0mu}}
\def\ripsqueeze{{\mskip -3.0mu}}
\def\ipcomma{\mathrel{,}}
\newbox\ipstrutbox
\setbox\ipstrutbox=\hbox{\vrule height8.5pt
width 0pt}
\def\ipstrut{\copy\ipstrutbox}
\def\lip#1<#2,#3>{\mathopen{\relax_{\ipstrut\ipscriptstyle{
#1}}\lipsqueeze
\langle} #2\ipcomma #3 \rangle}
\def\blip#1<#2,#3>{\mathopen{\relax_{\ipstrut
\ipscriptstyle{ #1}}\lipsqueeze\bigl\langle} #2\ipcomma #3 \bigr\rangle}
\def\rip#1<#2,#3>{\langle #2\ipcomma #3
\rangle_{\ripsqueeze\ipstrut\ipscriptstyle{#1}}}
\def\brip#1<#2,#3>{\bigl\langle #2\ipcomma #3
\bigr\rangle_{\ripsqueeze\ipstrut\ipscriptstyle{#1}}}
\def\angsqueeze{\mskip -6mu}
\def\smangsqueeze{\mskip -3.7mu}
\def\trip#1<#2,#3>{\langle\smangsqueeze\langle #2\ipcomma #3
\rangle\smangsqueeze\rangle_{\ripsqueeze\ipstrut\ipscriptstyle{#1}}}
\def\btrip#1<#2,#3>{\bigl\langle\angsqueeze\bigl\langle #2\ipcomma
#3
\bigr\rangle
\angsqueeze\bigr\rangle_{\ripsqueeze\ipstrut\ipscriptstyle{#1}}}
\def\tlip#1<#2,#3>{\mathopen{\relax_{\ipstrut\ipscriptstyle{
#1}}\lipsqueeze \langle\smangsqueeze\langle} #2\ipcomma #3
\rangle\smangsqueeze\rangle}
\def\btlip#1<#2,#3>{\mathopen{\relax_{\ipstrut\ipscriptstyle{
#1}}\lipsqueeze
\bigl\langle\angsqueeze\bigl\langle} #2\ipcomma #3
\bigr\rangle\angsqueeze\bigr\rangle}
\def\LUproperty{LU}
\def\lu{\ifmmode\LUpropery\else$\LUproperty$\fi}

\def\om{\omega}
\def\tg{\operatorname{tg}}
\def\RR{\mathbb{R}}
\def\CC{\mathbb{C}}
\def\TT{\mathbb{T}}
\def\ZZ{\mathbb{Z}}

\def\hA{\hat A}

\def\H{\mathcal{H}}

\def\ZM{\Zcenter \M}
\def\MyMath#1{\ifmmode #1\else$#1$\fi}

\def\M{M}
\def\UM{\mathcal{U}\M}
\def\ZM{\mathcal{Z}\M}
\def\ZUM{\mathcal{Z}\UM}
\def\cox{\MyMath{C_0(X)}}
\def\coxalg{\cox-algebra}
\DeclareMathOperator{\Glimm}{Glimm}
\def\gla{\ifmmode\Glimm(A)\else$\Glimm(A)$\fi}

\def\ewcr{\mathcal{CR}}
\def\CR{\MyMath{\ewcr}}

\def\LI{\MyMath{\mathcal{LI}}}

\def\tensor{\otimes}
\def\mtensor{\tensor_{\text{\normalfont max}}}
\def\atmb{A\mtensor B}

\def\xtensor{\tensor_X}
\def\coxk{C_0(X,\K)}

\def\E{\mathcal{E}}

\makeatletter
\def\labelenumi{(\@alph\c@enumi)}
\def\theenumi{\@alph \c@enumi}
\newcount\charno
\def\alphapart#1{\charno=96
\advance\charno by#1\char\charno}
\def\partref#1{\textup{(}\alphapart{#1}\textup{)}}
\makeatother

\def\hatu{\hat u}

\def\fc#1*#2{{#1} * {#2}}   
\def\guy{G\backslash Y}

\def\indgza{Z\times_G A}
\def\inda{\Ind\alpha}
\def\A{\mathcal{A}}
\def\lbb{\lbrack\!\lbrack}
\def\rbb{\rbrack\!\rbrack}
\def\eclass(#1,#2){\lbb #1,#2 \rbb}

\def\HP{\operatorname{HP}}

\def\dual(#1,#2){\spec{#1,#2}}

\DeclareMathOperator{\ad}{ad}

%% file: cxcab.tex
%
%

\begin{abstract}
Suppose that $G$ has a representation group $H$, that $\gab:=
G/\overline{[G,G]}$ is compactly generated, and that $A$ is a
\cs-algebra for which the complete regularization of $\Prim(A)$ is a
locally compact Hausdorff space $X$.  In a previous article, we showed
that there is a bijection $\alpha\mapsto(Z_\alpha,f_\alpha)$
between the collection of exterior equivalence classes of
locally inner actions
$\alpha:G\to\Aut(A)$, and the collection of principal $\hgab$-bundles
$Z_\alpha$ together with continuous functions $f_\alpha:X\to
H^2(G,\T)$.
In this paper, we compute the crossed products $A\rtimes_\alpha G$ in
terms of the data $Z_\alpha$, $f_\alpha$, and~$\cs(H)$.
\end{abstract}

%% file: intro1.tex
%
%
\section{Introduction}

This paper is a continuation of our study of locally inner actions
begun in \cite{ew2}.  In that article we gave a classification up to
exterior equivalence of actions of a smooth group $G$ on a
$\CR$-algebra $A$.  In this paper, we want to consider the structure
of the corresponding crossed products.

As in \cite{ew2}, we are motivated by a desire to make considerable
progress along the lines of a research program outlined by Rosenberg
in his
survey article \cite{ros5} (see ``Research Problem~1'' in \S3 of that
article).
As detailed there, it is interesting and important
to obtain information about crossed products
 of  actions with ``single orbit
type'' acting on  continuous-trace \cs-algebras.  Using the
Packer-Raeburn
stabilization trick, an action of $G$ on a continuous-trace \cs-algebra
$A$ with a single orbit type and constant stabilizer $N$ can be
decomposed into a spectrum fixing action of $N$ and an action of $G/N$
which acts freely on $\spec{A\rtimes_\alpha N}$.  Thus an important
first step in this program is to consider
spectrum fixing automorphism groups.  As it turns
out, provided
that the quotient
$\gab$ of
$G$ by the closure of its commutator subgroup $[G,G]$ is compactly
generated, spectrum fixing automorphism groups of continuous-trace
\cs-algebras $A$ are necessarily locally inner in that each point in
$\hA$
has a \nbhd{} $U$ such that the action restricts to an inner
action of the ideal $A_U$.  (This
follows from the proof of \cite[Corollary~2.2]{ros2}.)
Thus it is natural to try to classify locally
inner actions on arbitrary \cs-algebras rather than restricting ourselves
to actions on continuous-trace algebras.
In \cite{ew2}, for suitable $G$, we were able to do precisely this
for a large class of \cs-algebras: namely those
algebras whose primitive ideal space $\Prim(A)$ has a second countable
locally compact complete regularization $X$.
Following \cite[Definition~2.5]{ew2}, such algebras are called
\CR-algebras.
The collection of \CR-algebras whose primitive ideal spaces have
complete regularization (homeomorphic to) $X$ is denoted by $\CR(X)$.
As noted in \cite[\S2]{ew2}, all unital \cs-algebras are \CR-algebras,
as are all quasi-standard
algebras (as studied in \cite{as}).

In this paper, we want to give a precise bundle-theoretic description
of the crossed products corresponding to the dynamical systems
classified in \cite{ew2}.
As in \cite{ew2}, our methods require that (virtually) everything in
sight be separable.  Thus we assume from the onset that all our
automorphism groups are second countable, and that the \cs-algebras on
which they act be separable.  In particular, we shall make considerable
use of Moore's group cohomology groups $H^n(G,A)$ which are defined
for a second countable locally compact group $G$ and a Polish
$G$-module $A$ \cite{moore3,moore4}.  It is important to keep in
mind that the $H^n(G,A)$ are themselves topological groups --- at
least in the case $n=2$ and $A=\T$.
While the Moore topology on even $H^2(G,\T)$ can be pathological
without additional hypotheses on $G$, we assume, as we did in
\cite{ew2},
that $G$ is smooth which, among other things, implies that
$H^2(G,\T)$ is a locally compact abelian group.  A group is
smooth precisely when it has a representation group; although we will
give the definition of a representation group shortly, we point out
that the category of (second countable) smooth groups includes all
connected simply connected Lie groups, all compact groups, all
discrete groups, and all compactly generated abelian groups (see
\cite[Remark~4.2 and Corollary~4.6]{ew2}).

In Theorem~6.3 of \cite{ew2}, we showed that under the assumptions
that
$G$ is smooth, that
$\gab$ is compactly
generated, and that $A\in
\CR(X)$, then  the collection
$\LI_G(A)$ of exterior equivalence classes of
locally inner actions of $G$ on $A$ is parameterized by
\[
H^1(X,\shgab)\oplus C\(X,H^2(G,\T)\).
\]
If $(A,G,\alpha)$ is such a locally inner system,
our main result (\thmref{thm-main}) gives a description of the
crossed product $A\rtimes_\alpha G$ in terms of the associated
invariants $\zeta_H(\alpha)\in H^1(X,\shgab)$, $f_\alpha\in
C\(X,H^2(G,\T)\)$, and a representation group $H$ for $G$ as
described below.

The function $f_\alpha:X\to H^2(G,\T)$ arises naturally.  A
\CR-algebra is naturally a \coxalg, and therefore admits a fibering
over $X$ (see \secref{sec-bund-op}.1); $f_\alpha(x)$ is defined to be
the inverse of the Mackey obstruction for the induced action
$\alpha^x$ on the fibre $A_x$ (see Definition~\ref{def-zetaphi}).  The
construction of $\zeta_H(\alpha)$ is more subtle, although it reduces
to the usual Phillips-Raeburn obstruction (\cite{pr2},
\cite[\S2]{ew2}) when $\alpha$ is locally unitary.  As indicated by
the notation, it may depend on the choice of a representation group
for $G$.  A summary of the basic facts about smooth groups and
representation groups is given in
\cite[\S4]{ew2}.  We state some of the basic results here for
convenience.
Recall that if
\begin{equation}\label{eq-h}
e\arrow{e} C \arrow{e} H \arrow{e} G \arrow{e} e
\end{equation}
is a locally compact central extension of $G$ by an abelian group $C$,
then any Borel section $c:G\to H$ satisfying $c(eC)=e$ determines a
cocycle
\begin{equation}\label{eq-sigma1}
\sigma(s,t)=c(s)c(t)c(st)^{-1}
\end{equation}
in the Moore group $Z^2(G,C)$.  If $\chi\in\widehat C$, then
$\sigma_\chi = \chi\circ\sigma$ is a cocycle in $Z^2(G,\T)$.
The resulting map $\tg:\widehat C\to H^2(G,\T)$ is a continuous \hm{}
with respect to the Moore topology on $H^2(G,\T)$, and depends only on the
extension.  The map $\tg$ is called the \emph{transgression map}.
Moore called $G$ \emph{smooth} if $G$ has a central extension
\eqref{eq-h}, called a \emph{representation group},
 for which the transgression map is an isomorphism of
topological groups.  In that case, we can view $f_\alpha$ as a
continuous map from $X$ to $\widehat C$.  Since $C$ is central in $H$,
$\cs(H)$ admits a natural $C_0(\widehat C)$-action; that is, $\cs(H)$
is a $C_0(\widehat C)$-algebra.  The pull back
$f_{\alpha}^*\(\cs(H)\):= C_0(X)\tensor_{C_0(\widehat C)}\cs(H)$ is
then a \coxalg.
Since $A$ is also a \coxalg,
we can form the balanced tensor product
\begin{equation}\label{eq-mt}
A\tensor_f \cs(H):= A\tensor_{C_0(X)} f_\alpha^*\(\cs(H)\).
\end{equation}
In the special case that $\zeta_H(\alpha)$ is trivial, our main
theorem implies that $A\rtimes_\alpha G$ is isomorphic to
\eqref{eq-mt}.
When $\zeta_H(\alpha)$ is nontrivial, then it is necessary to
``twist'' \eqref{eq-mt} by a principal
$\hgab$-bundle $Z$ over $X$ whose isomorphism class corresponds to
$\zeta_H(\alpha)$ in $H^1(X,\shgab)$.  The details of this
construction are given in \secref{sec-bund-op.6} (see
Definition~\ref{def-product2}).  The basic idea is to view a \coxalg,
such as $A\tensor_f \cs(H)$, as a $\hgab$-bundle over $X$ and form what
corresponds to the usual bundle product: $\fc Z*{\(A\tensor_f\cs(H)\)}$.
The latter is naturally a \coxalg{} which admits a $\hgab$-action
which we denote $\fc
Z*\alpha$.  Then our main result goes as follows.

\begin{thmnn}
[{\thmref{thm-main}}]
Let $G$ be a smooth group with representation group $H$.  Suppose that
$\gab$ is compactly generated, that $A\in\CR(X)$, and that
$\alpha:G\to \Aut(A)$ is a locally inner action.
If $f_\alpha:X\to H^2(G,\T)$ and $\zeta_H(\alpha)$ are as above, and
if $q:Z\to X$ is a principal $\hgab$-bundle corresponding to
$\zeta_H(\alpha)$, then
there exists a $C_0(X)$-linear and
$\hgab$-equivariant isomorphism between
$A\rtimes_{\alpha}G$ and $\fc Z* {(A\otimes_{f_\alpha}C^*(H))}$.
\end{thmnn}

Another special case of interest arises when $\alpha$ is locally
unitary.  Then $f_\alpha$ is trivial and $\zeta_H(\alpha)$ is the
(generalized) Phillips-Raeburn obstruction.  (In fact, we do not
require $G$ to be smooth in this event.)  Then \thmref{thm-loc-uni}
implies that
\[
A\rtimes_\alpha G\cong A\tensor_{C_0(X)} \big(\fc Z*{\(C_0(X,\cs(G)\)}\big).
\]

Our basic motivation, and our basic strategy, for proving our results
involves viewing \coxalg s as the \cs-analogue of topological bundles
over $X$.  Thus we begin in \secref{sec-bund-op} with a review of some
of the basic facts about bundle operations and their \cs-counterparts:
\coxalg s, restrictions of
\coxalg s, balanced tensor products, and pull-backs.
In \secref{sec-bund-op.6}, we give our basic
product constructions alluded to above.

\secref{first-step} is devoted primarily to crossed products by inner
actions.  However, our methods require extensive use of the theory of
Busby-Smith twisted crossed products.  Much of the basis of this
theory has been worked out by Packer and Raeburn
\cite{para1,para2,para3}.  We review some of the basic facts here and then
formulate our results for twisted systems.

In \secref{loc-fibre}, we consider twisted systems which are ``locally
equivalent.''
Our main result here (\thmref{thm-mostimportant}) is crucial and
allows us to tie our analysis of inner systems to locally inner
systems.  As a rather special case, we derive the result on locally
unitary actions mentioned above (\thmref{thm-loc-uni}).   Our main
results on locally inner actions are spelled out in
\secref{loc-inner}.

In \secref{action}, we consider the special case of $\R^n$-actions.
Here the special structure of $\R^n$ allows us to give more detailed
information about the crossed products.
In a future article, we plan
to turn our attention to twisted transformation groups --- such as
arise in our study here (see \corref{cor-inner}).  This leads
naturally to the study of the group \cs-algebras of central extensions
in view of \lemref{lem-centralcross}(a).

%% file: bund-op1.tex
%
%

\let\paragraph\subsection

\section{Bundle operations on $C_0(X)$-algebras}\label{sec-bund-op}

It has long been known that many standard results in the theory of
\cs-algebras can
be motivated by viewing general \cs-algebras as function algebras over
``noncommutative spaces.''  This is especially true for the
constructions we require here.  In this section, we will review some
\cs-algebraic constructions --- most of which are well-known --- and
we want to emphasize that our constructions parallel the usual
topological notions of bundles, fibre products, pull-backs, and
$G$-fibre products.

\paragraph{\boldmath\coxalg s}
A (topological) \emph{bundle} $Y$ over a locally compact base space
$X$ is simply a continuous map $p:Y\to X$.  The fibres of the bundle
are $\set{p^{-1}(\{x\})}_{x\in X}$, and two bundles $p:Y\to
X$ and $ q:Z\to X$ over $X$ are {\em isomorphic}, if there exists a
homeomorphism $h:Y\to Z$ satisfying $q\circ h=p$.  This implies that
$h$ maps each fibre $p^{-1}(\{x\})$ homeomorphically onto
$q^{-1}(\{x\})$.

As we shall see,
the $\cs $-algebra analogue of a fibre bundle over $X$ is a
\emph{$C_0(X)$-algebra $A$}; that is,
 a $\cs $-algebra $A$ together with a
$*$-homomorphism $\phi$ from $C_0(X)$ to the center $\ZM(A)$
of the multiplier algebra $M(A)$ of $A$, which is nondegenerate in
that
\[
\phi\(\cox\)\cdot A:=\overline{\sp}\set{\phi(f)a:\text{$f\in\cox$ and $a\in
A$}}=A.
\]
Thus $A$ is a nondegenerate central Banach $\cox$-module
and we will usually suppress the map $\phi$ and
 write $f\cdot a$ in place of $\phi(f)a$.
If $A$ and $B$ are $C_0(X)$-algebras, then a homomorphism $\Psi:A\to B$
is called {\em $C_0(X)$-linear} if $\Psi(f\cdot a)=f\cdot \Psi(a)$ for
all $f\in C_0(X)$ and $a\in A$.
Two $C_0(X)$-algebras $A$ and $B$ are {\em isomorphic},
if there exists a $C_0(X)$-linear isomorphism $\Psi:A\to B$.

\coxalg s have enjoyed considerable interest of late, and there are
several nice treatments available \cite{may1,blanchard}.  We record
some of their basic properties here for convenience.
If $A$ is a \coxalg{}, if
$U$ is open in $X$, and if $J$ is the ideal of functions in \cox{}
vanishing off $U$, then
\begin{equation}
\label{eq-jdota}
J\cdot A:=\overline{\sp\set{f\cdot a:\text{$f\in J$ and $a\in A$}}}
\end{equation}
is an ideal in $A$.
The fibre $A_x$ of
$A$ over $x$ is defined to be the quotient
$A_x=A/I_x$, where $I_x:= C_0(X\setminus\{x\})\cdot A$.
The spectrum $\hA$
can then be written as a disjoint union $\coprod_{x\in X}\hA_x$,
and the projection $p:\hA\to X$ is a continuous map.
Conversely, if $p:\hA\to X$ is any continuous map,
and if we identify $C^b(\hA)$ with $\ZM(A)$ via the
Dauns-Hofmann Theorem, then $p$ induces a nondegenerate $*$-homomorphism
$\phi:C_0(X)\to C^b(\hA)\cong \ZM(A)$ by defining
$\phi(f)=f\circ p$, and then $p$ coincides with the projection corresponding
to the $C_0(X)$-structure on $A$ induced by $\phi$.

A $C_0(X)$-algebra can be viewed
as the algebra of sections of an (upper-semicontinuous)
$\cs $-bundle over $X$ as follows.
For each $x\in X$  and  $a\in A$,  let $a(x)$ denote the
image of $a$ in the fibre $A_x=A/I_x$. Then we have a faithful
representation of $A$ into the \cs-direct sum
$\bigoplus_{x\in X}A_x$ given by $a\mapsto
(a(x))_{x\in X}$. The set of sections $x\mapsto a(x)\in A_x$ for $a\in
A$,
satisfy
\begin{description}
\item[\normalfont (C-1)]
 $\|a\|=\sup_{x\in X}\|a(x)\|$;
\item[\normalfont (C-2)] $x\mapsto \|a(x)\|$ is upper semicontinuous and
vanishes at infinity --- that is,
\[
\set{x\in X:\|a(x)\|\ge \epsilon}\quad\text{is compact for all $\epsilon>0$;}
\]
\item[\normalfont (C-3)]
$(f\cdot a)(x)=f(x)a(x)$ for all $f\in C_0(X)$ and $a\in A$;
\item[\normalfont (C-4)]
 $\{a(x): a\in A\}=A_x$ for all $x\in X$.
\end{description}
Conversely, if $\set{A_x}_{x\in X}$ is a family of $\cs $-algebras (zero or
nonzero),
then any $\cs $-subalgebra of $\bigoplus_{x\in X}A_x$ which is closed under
pointwise multiplication with elements of $C_0(X)$, and which satisfies
Conditions (C-2)~and (C-4) above, becomes a $C_0(X)$-algebra by defining
the $C_0(X)$-action on $A$ by pointwise multiplication.
A $C_0(X)$-algebra is called a continuous $C_0(X)$-bundle if the
maps $x\mapsto \|a(x)\|$ are continuous for all $a\in A$.
By Lee's theorem \cite{lee2}, this is
equivalent to saying that the projection $p:\hA\to X$ is open.

\paragraph{Restrictions}
If $p:Y\to X$ is a topological bundle, and if $U$ is a subspace of
$X$, then the restriction of $p$ to $p^{-1}(U)$ defines a bundle over
$U$.  There is a nice analogue of this construction for \coxalg s
provided $U$ is locally compact.

\begin{definition}
\label{def-restr}
Let $A$ be a $C_0(X)$-algebra and let
$Y$ be a nonempty locally compact
subset of $X$. Let
\[C_0(Y)\cdot A:=\set{b\in \textstyle{\bigoplus_{y\in Y}}A_y:
\text{$b(y)=f(y)a(y)$
for some
$f\in C_0(Y)$ and $ a\in A$}}.
\]
Then $A_Y:=C_0(Y)\cdot A$
is called the {\em restriction} of $A$ to $Y$.
\end{definition}

While it is not immediately clear that $A_Y$ is even a subspace of
$\bigoplus_{y\in Y} A(y)$, the next result shows that $A_Y$
is a $C_0(Y)$-algebra.  The proof uses the Cohen
Factorization Theorem \cite{cohen} which implies that every element
$a$ in a nondegenerate
Banach $\cox$-module $B$ is of the form $f\cdot b$ for some
$f\in \cox$ and $b\in B$.
A nice proof of a version of Cohen's result
sufficient for our purposes can be found in
\cite[Proposition~1.8]{blanchard}.

\begin{lem}
\label{lem-exact}
Suppose that $Y$ is a nonempty locally compact subset of $X$, and that
$A$ is a \coxalg.  Then
the restriction $A_Y$ is a $C_0(Y)$-algebra with $(A_Y)_y=A_y$ for all
$y\in Y$.  If $U$ is open in $X$, then $A_U$ can be identified with
the ideal $C_0(U)\cdot A$ defined in \eqref{eq-jdota}.  If $C$ is
closed in $X$, then $A_C$ is the image of $A$ by the natural map of
$\bigoplus_{x \in X}A_x$ onto $\bigoplus_{x\in C} A_x$.  Moreover,
\begin{equation}
\label{eq-exactseq}
0\arrow{e} A_{X\setminus C}\arrow{e} A \arrow{e} A_C \arrow{e} 0
\end{equation}
is an exact sequence of \cs-algebras.
\end{lem}

\begin{proof}
We identify $C_0(U)$ with an ideal in \cox.  Note that $C_0(U)\cdot A$
is a nondegenerate Banach $C_0(U)$-module.  Therefore the Cohen
Theorem implies that every element of
$C_0(U)\cdot A$ is of the form $f\cdot a$ for $f\in C_0(U)$ and $a\in
A$.  Thus we can identify $A_U$ with $C_0(U)\cdot A$ as claimed.

The
assertion about closed sets follows from the Cohen Theorem applied to
$A$ together
with the observation that  restriction defines a
surjection  of $C_0(X)$ onto $C_0(C)$.

To establish \eqref{eq-exactseq}, it suffices to see that if $a\in A$
and if $a(x)=0$ for all $x\in C$, then $a\in C_0(X\setminus C)\cdot
A$.  But this is clear from the compactness of
\[
K=\set{x\in X:\|a(x)\|\ge\epsilon}
\]
for all $\epsilon>0$ so that $a$ can be approximated by elements of
the form $g\cdot a$ with $g\in C_c(X\setminus C)$.

It now remains only to prove the first assertion.  However, the first
part of the proof
shows that the first assertion is true if $Y$ is either open or
closed.  Since a subset of a locally compact space is locally compact
only if it is locally closed, the result follows.  In fact, $A_Y$ can
be identified with an ideal in $A_{\overline Y}$.
\end{proof}

The following lemma will be useful in the sequel.

\begin{lem}\label{lem-useful}
Let $A$ be a $C_0(X)$-algebra and let $B\subset \bigoplus_{x\in X}A_x$
be such that each $b\in B$ satisfies condition~{\normalfont (C-2)},
and such that for
each $x\in X$ and $b\in B$ there exists an open neighborhood $U$ of
$x$ such that $C_0(U)\cdot b=\set{f\cdot b:f\in C_0(U)}\subset A$.
Then $B\subseteq A$. If, in addition,
$B$ is a $C_0(X)$-submodule of $\bigoplus_{x\in X}A_x$
 such that for each $x\in X$ there
exists an open neighborhood $U$ of $x$ satisfying $C_0(U)\cdot B=A_U$,
then $B$ is a dense subspace of $A$.
\end{lem}
\begin{proof} Let $b\in B$ and let $g\in C_c(X)$.
Then, using condition~(C-2) and a partition of unity
for $\supp( g)$, it is not hard to show
that $g\cdot b\in A$.
Then, by taking an approximate unit in $C_0(X)$ which lies in
$C_c(X)$, we see that
condition~(C-1)  implies $b\in A$.

Assume that $B$ satisfies the additional assumptions.
Again using a partition of unity, it follows that each
$a\in A$ with $\supp( a)$  compact is a linear combination of elements in $B$.
Thus $B$  contains the set of
all elements $a\in A$ with compact support.
Therefore $B$ is dense in $A$.
\end{proof}

\paragraph{Balanced tensor products}
If $p:Y\to X$ and $q:Z\to X$ are topological bundles over $X$, then
the fibred product $Y\times_X Z=\set{(y,z)\in Y\times Z:p(y)=q(z)}$ is
naturally a bundle over $X$.  The projection $Y\times_XZ\to
X$ given by $ (y,z)\mapsto p(y)\; (=q(z))$ gives $Y\times_XZ$ a
canonical structure as a fibre bundle over $X$ with fibres
$p^{-1}(\{x\})\times q^{-1}(\{x\})$.
As will be clear from \lemref{lem-diag}, the $C_0(X)$-algebra analogue of
fibre products is a balanced tensor product $A\xtensor B$ defined
below.  In the case of nuclear \cs-algebras, one can
proceed along the lines of \cite[\S1]{rw}.  The general case has been
considered by several authors
\cite{kw,gior-mingo,blanchard,blanchard2}.
Blanchard's treatment is sufficient for our purposes, and as he
makes clear, it is most
convenient to work with maximal tensor products.  We record some of
the basic constructions here for the readers convenience.
(Note that
Blanchard considers only $X$ compact.)

If $A$ and $B$ are \cs-algebras, then there are natural injections
$i_A:M(A) \to M(\atmb)$ and $i_B:M(B)\to M(\atmb)$ such that
\begin{gather*}
i_A(c)(a\tensor b)=ca\tensor b,\qquad i_B(d)(a\tensor b)=a\tensor
db,\text{ and} \\
i_A(a)i_B(b)=i_B(b)i_A(a)=a\tensor b
\end{gather*}
for all $c\in M(A)$, $d\in M(B)$, $a\in A$, and $b\in B$ (e.g.,
\cite[Lemma~T.6.1]{wegge}).  In particular, if $A$ is a \coxalg{} and $B$
is a $C_0(Y)$-algebra, then $\atmb$ is a $C_0(X\times Y)$-algebra.
Moreover, just as is shown in \cite[Corollaire~3.16]{blanchard}, the
nice behavior of the maximal tensor product with respect to quotients
implies that the fibres of $\atmb$ are exactly what is expected, and
explains our preference for the maximal tensor product.

\begin{lem}[Blanchard]
\label{lem-maxtensor}
Suppose that $A$ is a $C_0(X)$-algebra and $B$ is a $C_0(Y)$-algebra.
Then, via the canonical composition of maps
$$C_0(X\times Y)\to \ZM(A)\otimes \ZM(B)\to \ZM(A\mtensor B),$$
$A\mtensor B$ is a $C_0(X\times Y)$-algebra with fibres
$(A\mtensor B)_{(x,y)}$ isomorphic to $A_x\mtensor A_y$.
Moreover, for any elementary tensor $a\otimes b$ we have
$(a\otimes b)(x,y)=a(x)\otimes a(y)\in A_x\mtensor A_y$.
\end{lem}

If $A$ and $B$ are both \coxalg s, then composition
with $i_A$ and $i_B$ gives $\atmb$ two \coxalg{} structures.
Since any quotient of a \coxalg{} is a \coxalg,
the two \coxalg{} structures will coincide on a given quotient
exactly when elementary
tensors of the form $f\cdot a\tensor b-a\tensor f\cdot b$ are mapped
to zero.

\begin{definition}
[{cf., \cite{blanchard}}]
\label{def-CXtensor}
Let $A$ and $B$ be two $C_0(X)$-algebras and let $I$ be the closed
ideal of $A\mtensor B$ generated by
\[
\set{a\cdot f\otimes b-a\otimes f\cdot b: \text{$a\in A$, $b\in B$,
$f\in C_0(X)$}}.
\]
Then $A\xtensor B:=(A\mtensor B)/I$ equipped with the
$C_0(X)$-action given on the images $a\xtensor  b$ of elementary
tensors $a\otimes b$ by
\[
f\cdot (a\xtensor  b)=  f\cdot a\xtensor  b
= a\xtensor  f\cdot b
\]
is called the (maximal)
{\em $C_0(X)$-balanced tensor product\/} of $A$ and $B$.
\end{definition}

\begin{remark}
When neither $A$ nor $B$ is nuclear, then it is possible to form other
balanced tensor products.  A detailed account may be found in
\cite{blanchard2}.  For example, it is observed there that
the second assertion of
\lemref{lem-diag} is false for Blanchard's minimal tensor
product $\tensor^m_{C(X)}$,
which is defined using the spatial tensor product.
\end{remark}

As in the above definition, we will denote the
image of an elementary tensor $a\otimes b$ in
$A\xtensor B$ by $a\xtensor b$.  Notice that $\atmb$ also has a
\coxalg{} structure arising from viewing $\cox$ as a the quotient of
$C_0(X\times X)$ by the ideal $C_\Delta$ of functions vanishing on the
diagonal $\Delta:=\set{(x,x):x\in X}$.  Our next result shows that
this structure also induces the given structure on $A\xtensor B$, and
that $\atmb$ coincides with Blanchard's $A\tensor^M_{C(X)} B$ when $X$
is compact.

\begin{lem}
[cf., {\cite[Proposition~2.2]{blanchard2}}]
 \label{lem-diag}
Let $A$ and $B$ be $C_0(X)$-algebras. Then
$A\xtensor B$ is $C_0(X)$-isomorphic to the restriction
$(A\mtensor B)_{\Delta }$ of $A\mtensor B$ to
$\Delta =\{(x,x):x\in X\}$,
where the $C_0(X)$-structure on $(A\mtensor B)_{\Delta }$
is defined via the canonical homeomorphism $x\mapsto (x,x)$
between $X$ and $\Delta $.

In particular, each fibre
$(A\xtensor B)_x$ is isomorphic to $A_x\mtensor B_x$
and the image of an elementary tensor $a\xtensor  b$ in $A\xtensor B$
is given by the section $x\mapsto a(x)\otimes b(x)\in A_x\mtensor B_x$.
\end{lem}
\begin{proof}
Since $\Delta $ is closed in $X\times X$, it follows from
\lemref{lem-useful}
that
$(A\mtensor B)_{\Delta }$ is the quotient of $A\mtensor B$
by the ideal $J=C_{\Delta}\cdot (A\mtensor B)$.
Since
$$(a\cdot f\otimes b-a\otimes f\cdot b)(x,x)=
f(x)\big(a(x)\otimes b(x)-a(x)\otimes b(x)\big)=0$$
for all $f\in C_0(X)$, $a\in A$, and $b\in B$, it follows that the balancing
ideal $I$ given in Definition~\ref{def-CXtensor} is contained in
$J$. Conversely, since $C_{\Delta}$ is
the closed ideal of $C_0(X\times X)$ generated by
$\{hf\otimes g- f\otimes hg: h,f,g\in C_0(X)\}$ and since
$(hf\otimes g- f\otimes hg)\cdot(a\otimes b)\in I$ for all elementary
tensors $a\otimes b$, it follows that the quotient map
$A\mtensor B\to A\xtensor B$ maps $C_{\Delta}\cdot(A\odot B)$ to $\{0\}$.
But this implies that $J\subseteq I$.
This proves the first assertion, and the second assertion now follows
from the first and  Lemma~\ref{lem-maxtensor}.
\end{proof}

If $A$ and $C$ are \cs-algebras, then a \hm{} $\Phi:A\to M(C)$ is
called nondegenerate if $\Phi(A)\cdot C$ is dense in $C$.
Recall that the maximal tensor product $A\mtensor B$
enjoys the following  universal property:
If $C$ is a $\cs $-algebra and $\Phi_A:A\to M(C)$ and $ \Phi_B:B\to M(C)$
are commuting nondegenerate $*$-homomorphisms,
then there exists
a unique nondegenerate homomorphism
$\Phi_A\otimes\Phi_B:A\mtensor B\to M(C)$ satisfying
$(\Phi_A\otimes\Phi_B)(a\otimes b)=\Phi_A(a)\Phi_B(b)$
for all elementary tensors $a\otimes b\in A\mtensor B$.
Conversely,
if $\Psi:A\mtensor B\to M(C)$ is a nondegenerate $*$-homomorphism,
then $\Psi=\Psi_A\otimes\Psi_B$, where
$\Psi_A=\Psi\circ i_A$ and $\Psi_B=\Psi\circ i_B$
(see \cite[Proposition~4.7]{taktext}). A similar result is true
for the maximal balanced tensor product:

\begin{prop}\label{prop-universal}
Let $A$ and $B$ be $C_0(X)$-algebras.
If $\Phi_A:A\to M(C)$ and $\Phi_B:B\to M(C)$ are commuting nondegenerate
\cox-homomorphisms,
then there exists a unique nondegenerate \cox-homomorphism
$\Phi_A\xtensor \Phi_B:A\xtensor B\to M(C)$
such that $(\Phi_A\xtensor \Phi_B)(a\xtensor b)=\Phi_A(a)\Phi_B(b)$
for all elementary tensors $a\xtensor b\in A\xtensor B$.

Conversely,  if
$\Psi:A\xtensor B\to M(C)$ is a nondegenerate $C_0(X)$-linear
homomorphism, then $\Psi=\Psi_A\xtensor \Psi_B$ for unique pair
$(\Psi_A,\Psi_B)$ of commuting nondegenerate \cox-homomorphisms.
\end{prop}
\begin{proof}
Let $\Phi_A$ and $\Phi_B$ be as in the proposition.
By the universal
property of $A\mtensor B$, there exists a homomorphism
$\Phi_A \otimes\Phi_B:A\mtensor B\to M(C)$ satisfying
$(\Phi_A\otimes\Phi_B)(a\otimes b)=\Phi_A(a)\Phi_B(b)$ for all
$a\in A$, $b\in B$.
Since $\Phi_A(a\cdot f)\Phi_B(b)=\Phi_A(a)\Phi_B(f\cdot b)$,
it follows  that $a\cdot f\otimes b-a\otimes f\cdot b$
lies in the kernel of $\Phi_A \otimes\Phi_B$ for all $f\in C_0(X)$,
$a\in A$, $b\in B$.
Thus $\Phi_A \otimes\Phi_B$ factors through a nondegenerate
homomorphism $\Phi_A\xtensor \Phi_B:A\xtensor B\to M(C)$, which
satisfies $\Phi_A\otimes_X\Phi_B(a\xtensor
b)=\Phi_A\otimes\Phi_B(a\otimes b)=
\Phi_A(a)\Phi_B(b)$ for all elementary tensors $a\xtensor b$.

For the converse put $\Psi_A=\Psi\circ i_A$ and
$\Psi_B=\Psi\circ i_B$.
\end{proof}

\begin{remark}\label{rem-universal}
Suppose that $C$ is a \cs-algebra and that $A$ and $B$ are \coxalg s.
The above argument shows that if $\Phi_A:A\to M(C)$ and $ \Phi_B:B\to M(C)$
are nondegenerate homomorphisms such that $\Psi_A(f\cdot
a)\Psi(b)=\Psi_A(a)\Psi_B( f\cdot b)$ for all $a\in A$ and $b\in B$,
then there is a unique nondegenerate \hm{} $\Psi:A\xtensor B\to M(C)$
as in \propref{prop-universal}.
Conversely, if $\Phi=A\xtensor B$ is a nondegenerate \hm{} into
$M(C)$, then
$\Psi=\Psi_A\xtensor \Psi_B$ for a pair of commuting nondegenerate \hm
s such that $\Phi_A(a\cdot f)\Phi_B(b)=\Phi_A(a)\Phi_B(f\cdot b)$.
\end{remark}

\paragraph{Pull-backs}
If $p:Y\to X$ is a topological bundle and if $f:Z\to X$ is a
continuous map, then the pull-back $f^*Y=\set{(z,y)\in Z\times Y:
f(z)=p(y)}$ is the bundle $f^*p:f^*Y\to Z$, where $f^*p(z,y)=z$.  Note
that the fibre over $z$ can be identified with $p^{-1}\(f(z)\)$.  The
analogous construction for \cs-algebras turns out to be a certain
balanced tensor product and was introduced in \cite{rw}.
Here we confine ourselves to the following remarks.
\begin{remark}
\label{rem-final} (a)
If $A$ and $B$ are nuclear (and separable) \coxalg s, then
the balanced tensor product $A\xtensor B$ coincides with the
construction given by Iain Raeburn and the second author  in
\cite{rw}. In particular, it follows that
if $p:\Prim(A)\to X$ and $ q:\Prim(B)\to X$ are the projections
determined
by the $C_0(X)$-structures of $A$ and $B$, then
$\Prim(A\xtensor B)$ is homeomorphic to the fibre product
$\Prim(A)\times_X\Prim(B)$ \cite[Lemma~1.1]{rw}.
If $A$ or $B$ is type~I, then
we also have $\spec{A\xtensor B}\cong \hA\times_X\hat B$.

(b)
An important special case occurs when $B=C_0(Y)$
for a locally compact space $Y$. If $p:Y\to X$ is a continuous map,
(so that $C_0(Y)$ becomes a $C_0(X)$-algebra via the homomorphism
$\phi:C_0(X)\to C^b(Y)$ defined by $\phi(g)=g\circ p$), then
$A\xtensor C_0(Y)$ is not only a $C_0(X)$-algebra, but there is
also a canonical $C_0(Y)$-action on the balanced tensor product
given by the canonical embedding of $C_0(Y)$ into
$M(A\xtensor C_0(Y))$. We will write $A\otimes_pC_0(Y)$
for the balanced tensor product $p^*A:=A\xtensor C_0(Y)$ viewed
as a $C_0(Y)$-algebra; this is the {\em pull-back\/} of $A$ along $p$
as
defined in \cite{rw}.
If $y\in Y$, then the fibre $(A\otimes_pC_0(Y))_y$ is equal to $A_{p(y)}$,
and the projection $q:\spec{A\otimes_pC_0(Y)}\cong Y\times_X\hat{B}\to Y$
is given by $q(y,\pi)=p(y)$.
The justification for the pull-back terminology is
\cite[Proposition~1.3]{rw}
where is is shown that if $A$ is the
section algebra of a \cs-bundle, then $p^*A$ is the section algebra of
the pull-back bundle.

(c)
More generally, suppose that $B$ is a \cox-algebra and that
$A$ is a $C_0(Y)$ algebra.  Then if $p:Y\to X$ is continuous, we can
view $B$ as a $C_0(Y)$-algebra via composition with $p$.  Since
\begin{equation*}
B\tensor_Y A \cong \(B\tensor_p C_0(Y)\)\tensor_Y A,
\end{equation*}
we will write $B\tensor_p A$ in place of $B\tensor_Y A$.
\end{remark}

%% file: bund-op2.tex
%
%



\section{The bundle product constructions}
\label{sec-bund-op.6}

A
{\em topological bundle with group $G$}
is a topological bundle $p:Y\to X$ such that $Y$ is a $G$-space
in such a way that each $s\in G$ acts as a bundle
isomorphism of $p:Y\to X$.
The \cs-algebraic analogue of a bundle with group $G$, is a
\cs-dynamical system $(A,G,\alpha)$ in which $A$ is a \coxalg, and
each $\alpha_s$ is a \cox-automorphism.  More simply, $\alpha$ is a
strongly continuous
homomorphism of $G$ into $\Aut_{C_0(X)}(A)$.

As mentioned in the introduction, we want to build a dynamical system
$(\fc Z*A,G,\fc Z*\alpha)$ from a ``product'' of a \cox-system
$(A,G,\alpha)$ and a principal $G$-space $Z$.  Our construction is
based on the standard product of two fibre bundles with structure
group $G$.  Since the theory is considerably easier to describe in the
case when $G$ is abelian, and since the applications we require all
involve the groups $\hgab$, we will assume that $G$ is
abelian for the remainder of this section.

Recall that a {\em fibre bundle\/} is a topological
$G$-bundle in which each
of the fibres has been identified with a fixed topological space
$F$. If $p:Y\to X$ is a fibre bundle, then we say that
$p:Y\to X$ is a {\em fibre bundle with group $G$\/} if
$F$ is a $G$-space and
$p:Y\to X$ is a topological
$G$-bundle such that the induced action of $G$ on  each fibre
of $Y$ coincides with the given $G$-action on $F$ under the
identification of the fibres with $F$.
Such a bundle is {\em locally trivial\/}
if each point in $X$ has a
\nbhd{} $U$ such that $p^{-1}(U)$ is $G$-homeomorphic to $U\times F$
(where $G$ acts on the second factor).  If $F$ equals $G$ with $G$
acting by translation, then we
obtain a \emph{principal $G$-bundle}.  Notice that $X$ then has an
open cover $\set{U_i}_{i\in I}$ such that there are
$G$-homeomorphisms $h_i:U_i\times G\to p^{-1}(U_i)$ for each $i\in
I$.  On overlaps
$U_{ij}:=U_i \cap U_j$ we obtain continuous \emph{transition
functions} $\gamma_{ij}:U_{ij}\to G$ such that
\begin{equation}\label{eq-transfcns}
h_j^{-1}\circ h_i(x,s)=\(x,s\gamma_{ij}(x)\)\quad\text{for $x\in
U_{ij}$ and $s \in G$.}
\end{equation}
Then, if $x$ belongs to a triple overlap $U_{ijk}:=U_i\cap U_j\cap
U_k$,
\[
\gamma_{ij}(x)\gamma_{jk}(x)=\gamma_{ik}(x).
\]
Therefore $\gamma:=\set{\gamma_{ij}}_{i,j\in I}$ defines a $1$-cocycle
in $Z^1(X,\sheaf G)$, and it is well know
that the class of $\gamma$ in the sheaf cohomology group
$H^1(X,\sheaf G)$ depends only on the isomorphism class of $p:Y\to X$.
(If
$G$ is a topological group, we use the caligraphic
letter
$\sheaf{G}$ to denote the corresponding sheaf of germs of
$G$-valued functions.)
Furthermore, $H^1(X,\sheaf G)$ can be identified with the collection of
isomorphism classes of principal $G$-bundles over $X$ (cf., e.g.,
\cite[\S5.33]{warner} or \cite{steenrod}).
Since $G$ is abelian, $H^1(X,\sheaf G)$ is a group under
pointwise product and the bundle product we want to investigate is
analogous to the
construction on the principal $G$-bundles corresponding the the
product in $H^1(X,\sheaf G)$.

\begin{remark}
\label{rem-trans-sections}
There is a certain arbitrariness in our definition of the transition
functions $\set{\gamma_{ij}}$ classifying a principal bundle --- this
is especially true when working with abelian groups.  We have made an
effort to be consistent with the convention given in
\eqref{eq-transfcns}.
Thus if we have local sections $\phi_i:U_i\to Y$ to
the bundle map $p:Y\to X$, then we can define $h_i(x,s)=s\cdot
\varphi_i(x)$, and our conventions force
\begin{equation}
\label{eq-ddag}
\varphi_i(x)=\gamma_{ij}(x)\cdot \varphi_j(x)\quad\text{for all $x\in
U_{ij}$.}
\end{equation}
\end{remark}

If $G$ and $Y$ are locally compact, then it is not hard to see that if
$p:Y\to X$ is a principal $G$-bundle, then the action of $G$ on $Y$ is
free and proper; that is, $s\cdot y=y$ implies $s=e$, and the map
$(s,y)\mapsto (s\cdot y,y)$ is proper as a map from $G\times Y$ to
$Y\times Y$.  Moreover, $p$ induces a homeomorphism of $\guy$ onto
$X$.  On the other hand, if $Y$ is a free and proper $G$-space and if
$p:Y\to X$ identifies $\guy$ with $X$, then $p:Y\to X$ is called a
\emph{proper $G$-bundle} \cite[\S2]{90a}.  A proper
$G$-bundle is a principal $G$-bundle exactly when there are local
(continuous) sections for $p:Y\to X$ \cite[Proposition~4.3(3)]{rw2}.
If $G$ is a
Lie group, then it follows from Palais's slice theorem \cite{palais}
that the proper $G$-bundles are precisely the principal $G$-bundles.
Of course, there are groups $G$ for which there exist non-principal proper
$G$-bundles
\cite[Remark~2.5]{90a}.

\begin{definition}
\label{def-product}
 Let $G$ be an \emph{abelian}  locally compact group and let
$q:Z\to X$ be a proper $G$-bundle over $X$.
\begin{enumerate}
\item  If $Y$ is any $G$-space,
then we define $Z\times_G Y$ to be the orbit space
$G\backslash (Z\times Y)$, where the $G$-action is defined by
$s\cdot (z,y):= (sz, s^{-1}y)$. We define a continuous map
$i:Z\times_G Y\to X$ by $i([z,y]):=q(z)$.  We define a left $G$-action
on  $Z\times_GY$ by $s\cdot [x,y]:=[s\cdot x,y]$, where
$[z,y]$ denotes the orbit of $(z,y)\in Z\times Y$.
\item  If $p:Y\to X$
is any topological bundle over $X$ with group $G$,
we define $r:Z*Y\to X$ to be the topological bundle
over $X$ with group $G$ such that $Z*Y:=\{[z,y]\in
Z\times_GY:q(z)=p(y)\}$ and $r:= i\restr{Z*Y}$.   The $G$-action is
induced from the $G$-action on $Z\times_GY$. We call
$r:Z*Y\to X$ the {\em $G$-fibre product\/} of $Z$ and $Y$.
\end{enumerate}
\end{definition}

\begin{remark}
\label{rem-product1}
Notice that the above definitions of $Z\times_GY$ and $Z*Y$ only make
sense when $G$ is abelian. For non-abelian $G$ one could define
similar spaces by taking the quotient of $Z\times Y$ by the diagonal
action. However, there would be no analogue of the $G$-actions on
$Z\times_GY$ and $Z*Y$.
\begin{enumerate}
\item
It is straightforward to check that $Z\times_GY$  is
a fibre bundle over $X$ with group $G$ and fibre $Y$ and that
$ Z*Y$ is a topological bundle with group $G$. The
$G$-isomorphism classes of $Z\times_GY$ and $Z*Y$ depend only on
those of $Z$ and $Y$.
\item If $Y$ is topological bundle over $X$ with group $G$ and
if $z\in Z$, then $[z,y]\mapsto y$ defines a $G$-equivariant isomorphism
of $r^{-1}\(\set{r([z,y])})$ onto $p^{-1}\(\set{r([z,y])})$.
Thus, $Z*Y$ has the same fibres as does $Y$, and we can view $ Z*Y$
as the bundle $Y$ \emph{twisted} by $Z$.
\item
If   $Z$ and $Y$ are proper
$G$-bundles, then it was shown in \cite[Lemma~2.4]{90a} that $Z*Y$ is
a proper $G$-bundle and that
$[Z][Y]\mapsto [Z*Y]$ defines an abelian group structure on the set
$\HP(X,\mathcal G)$ of all isomorphism classes of proper $G$-bundles
over $X$. The identity is given by the class of the trivial bundle
$X\times G$ and  the inverse of (the class of) $q:Z\to X$ is given by
(the class of) $\bar{q}:\bar{Z}\to X$, where $\bar Z:=\{\bar z: z\in
Z\}$ is identified with $Z$ as a topological space and
$s\cdot\bar z= (s^{-1}\cdot z)^{-}$. If $q:Z\to X$ and
$p:Y\to X$ are  principal bundles corresponding to
the classes $[\gamma_1]$ and $ [\gamma_2]$ in $H^1(X,\mathcal G)$, then
$r:Z*Y\to X$ is a principal $G$-bundle corresponding to
$[\gamma_1\cdot \gamma_2]$. Thus $H^1(X,\mathcal G)$ can be viewed
as a subgroup of
$\HP(X,\mathcal G)$ \cite[Remark~2.7]{90c}.
\item The general construction $Z*Y$ can be used for to obtain the
classification of locally trivial bundles with structure group $G$ by
$H^1(X,\mathcal G)$ via the classification of the principal $G$-bundles:
If $G$ acts effectively on the constant fibre $F$, and if $q:Z\to X$ is
a principal $G$-bundle corresponding to $[\gamma]
\in H^1(X,\mathcal G)$,
then $r:Z*Y\to X$ is a locally trivial bundle with fibre $F$
corresponding to $[\gamma]$, where $Y$ is the
trivial bundle $X\times F$. Since we will not require this
construction here, we omit the proof.
Note that
$[z, f]\mapsto [z,
(q(z),f)]$ is always a $G$-homeomorphism between
$Z\times_GF$ and $Z*Y$.
\end{enumerate}
\end{remark}

Now we introduce the \cs-algebraic constructions which will turn out
to parallel those
of $Z\times_G Y$ and $Z* Y$ above when $Y$ is replaced by a \coxalg{}
(see Propositions
\ref{prop-dualproduct}~and \ref{prop-bundle}).
Suppose that $q:Z\to X$ is a proper $G$-bundle and that
$(A,G,\alpha)$ is a \cs-dynamical system.  Notice that since $G$ is
abelian, $s\mapsto \alpha^{-1}_s$ is also a \hm.
Therefore we can employ the well-known construction of induced systems
\cite{rr,rw,ra2} to $(A,G,\alpha^{-1})$.
In particular, we define $\indgza$ to be $\Ind_G^Z(A,\alpha^{-1})$
(in the notation of
\cite{ra2}); that is, $\indgza$ is the
set of all bounded continuous functions $F:Z\to A$ satisfying
\[
\alpha_s(F(z))=F(s^{-1}\cdot z), \quad \text{for all $s\in G$ and $ z\in
Z$,}
\]
and such that $z\mapsto\|F(z)\|$ vanishes at
infinity on $X=G\backslash Z$. Equipped with the pointwise operations
and the supremum
norm, $\indgza$ becomes a $\cs$-algebra.
We define a strongly continuous action
$\Ind \alpha$ of $G$ on $\indgza$ by
\[
(\Ind\alpha)_s(F)(z)=\alpha_s(F(z))=F(s^{-1}\cdot z).
\]
Strong continuity follows from
straightforward compactness arguments
 using the fact that $z\mapsto
\|f(z)\|$ vanishes at infinity on $G\backslash Z$.

Note that $C_0(X)$ acts on $\indgza$ via
$(g\cdot F)(z)=g(q(z))F(z)$, $g\in C_0(X)$, so that
$(\indgza, G,\Ind\alpha)$ is actually a $C_0(X)$-system.
Note that the fibres $(\indgza)_x$ are all identified with $A$.

If $(A,G,\alpha)$ is itself a $C_0(X)$-system, then
there is a $C_0(X\times X)$-action on $\indgza$
given by
\[
(h\cdot F)(z)(x)=h(q(z),x)F(z)(x),\quad h\in C_0(X\times X).
\]
We define $\fc Z*A$ to be the restriction of
$\indgza$ to the diagonal $\Delta$ of $X\times X$.
Identifying $X$ with $\Delta$
gives $\fc Z*A$ the structure of a $C_0(X)$-algebra with fibres $(\fc
Z*A)_x\cong A_x$.
\begin{remark}
\label{rem-new}
\lemref{lem-exact} implies that
$\fc Z*A$ may be written as the
set of sections
\[
\set{f:Z\to\textstyle{\bigoplus_{x\in X}}A_x:\text{$
f(z)=F(z)\(q(z)\)\in
A_{q(z)}$ for some
$F\in \indgza$}}.
\]
Since $\inda$ is $C_0(X\times X)$-linear, it follows
that it restricts to a $C_0(X)$-linear action
$\fc Z*\alpha$ of $G$ on $\fc Z*A$.
If $F\in \indgza$, and $f(z)=F(z)\(q(z)\)$, then
\[
(\fc Z*\alpha)_s(f)(z)=\alpha_s\(F(z)\)\(q(z)\)=\alpha^{q(z)}_s\(f(z)\),
\]
where $\alpha^{q(z)}$ is the induced action on the fibre $A_{q(z)}$.
\end{remark}

\begin{definition}
\label{def-product2}
Suppose that $G$ is an abelian group and $q:Z\to X$ is a proper
$G$-bundle over $X$.
\begin{enumerate}
\item If $(A,G,\alpha)$ is a $\cs$-dynamical system then
$(\indgza, G,\inda)$ is called the $C_0(X)$-system {\em induced
from $(A,G,\alpha)$ via $Z$}.
\item If $(A,G,\alpha)$ is a
$C_0(X)$-system, then $(\fc Z*A,G,\fc Z*\alpha)$ is called the {\em
$G$-fibre product\/} of
$q:Z\to X$ with $(A,G,\alpha)$.
\end{enumerate}
\end{definition}

\begin{remark}\label{rem-product2}
(a)
If $\gamma$ is the diagonal action on $C_0(Z,A)$ given by
$\gamma_s(F)(z):= \alpha^{-1}_{s}\(F(s^{-1}\cdot z)\)$, then
$\indgza$ was denoted by $GC(Z,A)^\gamma$ in \cite{rw}
and is a generalized fixed point algebra for $\gamma$
\cite[Example~2.6]{rieffpr}.
The algebra $\fc Z*A$ defined above
first appeared in \cite[\S 2]{rw} as
$GC(Z,A)^{\gamma}/I$, where
$I$ is
the kernel of the quotient map $\indgza\to \fc Z*A$.
What is new with our construction are the actions $\inda$ and
$\fc Z*\alpha$.

(b) Our definitions of $\indgza$ and $A^Z$ only make
sense if $G$ is abelian, since we need $\alpha^{-1}$ to be
an action of $G$ on $A$.
Of course, we could have defined $\indgza$
as the algebra $\Ind_G^Z(A,\alpha)$.  This
would have the advantage of working for nonabelian $G$,
and would lead
to a sensible definition of $\fc Z*A$ in the general case.
However,
there would be no analogues for the actions
$\inda$ and $\fc Z*\alpha$ in the nonabelian case.
In any event, our definition more closely parallels the classical
bundle product, and leads to more elegant statements of our main
results.
\end{remark}

We now turn to the basic properties of the
$C_0(X)$-systems $(\indgza, G,\inda)$ and $(\fc Z*A, G, \fc Z*\alpha)$.
In so doing, we will see that these \cs-constructions from
Definition~\ref{def-product2} do
parallel the topological constructions
of Definition~\ref{def-product}.
The first result shows that
the construction of induced systems are
a special
case of the $G$-fibre product construction; therefore
we can focus on the latter.

\begin{prop}\label{proptriv}
Let $(A,G,\alpha)$ be any $\cs$-dynamical system with $G$ abelian,
and let $q:Z\to X$ be a proper $G$-bundle. Then
$(\fc Z*{C_0(X,A)},G,\fc Z*{(\id\otimes\alpha)})$ is canonically isomorphic
to $(\indgza, G, \Ind\alpha)$.
\end{prop}
\begin{proof}
Define $\Psi:Z\times_G C_0(X,A)\to \indgza$ by $\Psi(F)(z)=F(z,q(z))$.
Then it is easy to
check that
$\Psi$ is a surjective $*$-homomorphism, and Remark~\ref{rem-new}
implies that $\Psi$ induces an equivariant isomorphism between
$\fc Z*{C_0(X,A)}$ and $\indgza$.
\end{proof}

If $(A,G,\alpha)$ is a $C_0(X)$-system, then
$\hA$ is a topological bundle over $X$ with group $G$ with
respect to the projection $p:\hA\to X$ and the action
of $G$ defined by $s\cdot \pi=\pi\circ \alpha^{-1}_{s}$.
Now if $(z,\pi)\in Z\times \hA$, then it was shown in
\cite[Proposition 3.1]{rw} that $(z,\pi)$ determines an
irreducible representation $M(z,\pi)\in \specnp{(\indgza)}$
defined by $M(z,\pi)(F)= \pi(F(z))$. Moreover
$M(z_1,\pi_1)$ is equivalent to $M(z_2,\pi_2)$ if and only if
there exist an $s\in G$ such that $z_2= s\cdot z_1$ and
$\pi_2=s^{-1}\cdot \pi_1$ (note that $s^{-1}$ appears in the latter
formula as we have replaced $\alpha$ by $\alpha^{-1}$ in the formulae
from \cite{rw}).
The representations of
$\fc Z*A$ are then given by those $M(z,\pi)$ which satisfy
$q(z)=p(\pi)$. Thus we obtain

\begin{prop}[{\cite[Proposition~3.1]{rw}}]
\label{prop-dualproduct}
Let $q:Z\to X$ be a proper $G$-bundle with $G$ abelian.
\begin{enumerate}
\item
If $(A,G,\alpha)$ is a system, then
$\spec{\indgza}$ is naturally isomorphic to $Z\times_G\hA$ as
a fibre bundle over $X$ with group $G$.
\item
If $ (A,G,\alpha)$ is a $C_0(X)$-system, then
$\spec{\fc Z*A}$ is naturally isomorphic to $Z*\hA$ as topological
bundles over $X$ with group $G$.
\end{enumerate}
\end{prop}

Since $\indgza$ and $\fc Z*A$
are commutative if $A$ is commutative, the Gelfand theory yields an
immediate corollary:

\begin{cor}\label{cor-comm}
Let $Y$ be a locally compact $G$-space, and define $\tau:G\to
\Aut(C_0(Y))$ by $\tau_s(f)(y)=f(s^{-1}y)$.
Then
$Z\times_G C_0(Y)$ is equivariantly isomorphic to $C_0(Z\times_GY)$.
Moreover, if
$p:Y\to X$ is a locally compact topological bundle over
$X$ with group $G$, Then $\fc Z*C_0(Y)$ is $G$-equivariantly
isomorphic to $C_0(Z*Y)$.
\end{cor}

A similar result holds when $A$ is a continuous
\cox-bundle.

\begin{prop}[{cf., \cite[Proposition~2.15]{kmrw}}]\label{prop-bundle}
Suppose that $(A,G,\alpha)$ is a $C_0(X)$-system such that
$A$ is actually a continuous \cox-bundle; that is,
$A$ is the section algebra
$C_0(X;\A)$ of a \cs-bundle $p:\A\to X$. Then
$\A$ is a (continuous) $G$-space with the action
characterized by $s\cdot a(x) = \alpha_s(a)(x) = \alpha_s^x\(a(x)\)$,
$Z*\A$ is a $\cs$-bundle over $X$, and $\fc Z*A$ is canonically
isomorphic to $C_0(X, Z*\A)$.
\end{prop}
\begin{proof}
We omit the proof that
$\A$ is a continuous $G$-space with respect to the above
given action and that $\fc Z*\A$ is a $\cs$-bundle over $X$
(for more details see \cite[Proposition~2.15]{kmrw}).
In order to see that $\fc Z*A$ is isomorphic to $C_0(X,Z*\A)$
let $F\in \indgza$. Then $F$ defines a section $f_F\in C_0(X;
Z*\A)$ by
\[
f_F\(q(z)\)=[z,F(z)\(q(z)\)].
\]
The collection $\Gamma=\set{f_F:F\in\indgza}$ is dense in $C_0(X;
Z*\A)$ by \cite[Corollary~II.14.7]{fell-doran}, and it follows from
Remark~\ref{rem-new} that $\fc Z*A \cong C_0(X; Z*\A)$.
\end{proof}

We start to investigate the structure of $(\fc Z*A,G,\fc Z*\alpha)$
with some interesting special cases.  First notice that  if
$\alpha:G\to \Aut(A)$ is the trivial action, then
$(\fc Z*A,G,\fc Z*\alpha)$ is isomorphic to $(A,G,\alpha)$ for all
proper $G$-spaces $q:Z\to X$ (for a proof see
\cite[Proposition~3.2]{rw}).
A similar result holds when
$Z$ is a trivial bundle:

\begin{lem}
[{cf., \cite[Proposition 3.2]{rw}}]
\label{lem-trivial}
Let $(A,G,\alpha)$ be a $C_0(X)$-system
and let $q:Z\to X$ be a trivial $G$-bundle.
Let $\varphi:X\to Z$ be a continuous section for
$q:Z\to X$, and let $s(z)$ be the unique element in $G$
which satisfies $z=s(z)\varphi(q(z))$ for each $z\in Z$.
Then
$\Phi(f)(x)=f(\varphi(x))$ defines an equivariant
\cox-isomorphism of the systems
$(\fc Z*A,G,\fc Z*\alpha)$ and $(A,G,\alpha)$,
with inverse given by
$
\Phi^{-1}(a)(z)=\alpha^x_{s(z)^{-1}}\(a(q(z))\)$.
\end{lem}
\begin{proof}
Define $\Psi:\indgza \to C_0(X,A)$ by
$\Psi(F)(x)=F(\varphi(x))$.
Then it is easy to check that $\Psi$ is a $C_0(X\times X)$-linear
isomorphism with inverse given by
$\Psi^{-1}(g)(z)=\alpha^{-1}_{s(z)}\big(g(q(z))\big)$.
If $s\in G$ and $F\in
\indgza $, then
\[\Psi\((\inda) _s(F)\)(x)=(\inda) _s(F)(\varphi(x))=
\alpha_s\(F\(\varphi(x)\)\)
=\alpha_s\(\Phi(F)(x)\).
\]
Thus $\Psi$ carries $\inda $ to $\id\otimes\alpha$.
Since $\Psi$ is $C_0(X\times X)$-linear and the restriction
$(C_0(X,A)_{\Delta},G,(\id\otimes \alpha)^{\Delta})$ of
$(C_0(X,A), G, \id\otimes \alpha)$ to the diagonal $\Delta$
is isomorphic to $(A,G,\alpha)$, it follows that $\Psi$ induces
a $G$-equivariant and $C_0(X)$-linear isomorphism
$\Phi: \fc Z*A\to A$.
Evaluation at the fibres reveals that $\Phi$ and $\Phi^{-1}$
are given by the formulas in the lemma.
\end{proof}

\begin{remark}\label{rem-trivial}
The isomorphism of $\fc Z*A$ and $A$ given in
\lemref{lem-trivial} depends on the choice of section.
If $\Phi_1$ and $\Phi_2$ are induced from two different continuous sections
$\varphi_1$ and $\varphi_2:X\to Z$, then as in
\eqref{eq-ddag},
let $\gamma_{12}:X\to G$
be the transition function defined by
\[
\varphi_1(x)=\gamma_{12}(x)\varphi_2(x).
\]
Then
\begin{align*}
\Phi_2(f)(x)&= f\(\varphi_2(x)\) = F\(\varphi_2(x)\)(x) =
F\(\gamma_{12}(x)^{-1} \varphi_2(x)\)(x) \\
&= \alpha_{\gamma_{12}(x)}^x\(F\(\varphi_1(x)\)(x) =
\alpha_{\gamma_{12}(x)}^x \(f\(\varphi_1(x)\)\) \\
&= \alpha_{\gamma_{12}(x)}^x \(\Phi_1(f)(x)\).
\end{align*}
for all $x\in X$.
\end{remark}

Again, suppose that
$q:Z\to X$ is a proper $G$-bundle over $X$.
If $W$ is any locally compact subset of $X$, then the restriction
$Z_W:=q^{-1}(W)$ is a proper $G$-bundle over $W$, and our next result
shows that $\fc Z*A$ behaves well
with respect to restrictions.

\begin{lem}\label{lem-restr}
Let $(A,G,\alpha)$ be a $C_0(X)$-system and let $q:Z\to X$
be a proper $G$-bundle over $X$. If $W$ is a locally compact subset
of $X$, then $((\fc Z*A)_W,G,(\fc Z*\alpha)_W)$ and
$(\fc {Z_W}*{(A_W)}, G, \fc {Z_W}*{(\alpha^W)})$ are
isomorphic as $C_0(W)$-systems.
In particular, $((\fc Z*A)_x,G,(\fc Z*\alpha)^x)$
is isomorphic to $(A_x,G,\alpha^x)$
for all $x\in X$.
\end{lem}
\begin{proof}
The second assertion is a consequence of the first and
\lemref{lem-trivial}.
The first assertion is straightforward when $W$ is open
or closed. Since $Y$ is always the intersection of an open and a closed
set, the result follows by iteration.
\end{proof}

When
$Z$ is a principal $G$-bundle, it will be convenient to have
a description of $\fc Z*A$ in terms of
a representative $\gamma\in Z^1(X,\mathcal G)$ for the class in
$H^1(X, \mathcal G)$ corresponding to~$Z$.

\begin{prop}\label{prop-principal}
Let $(A,G,\alpha)$ be a $C_0(X)$-system with $G$ abelian and
$X$ paracompact. Let
$q:Z\to X$ be a principal bundle and let $\set{U_i}_{i\in I}$ be a locally
finite cover of $X$ such that $\gamma=\set{\gamma_{ij}}_{i,j\in I}$
represents the class in $H^1(X,\mathcal G)$ corresponding to $Z$.
Then a $C_0(X)$-system $(B,G,\beta)$ is \cox-isomorphic to
$(\fc Z*A,G,\fc Z*\alpha)$ if and only if there exist isomorphisms
$\Phi_i:B_{U_i}\to A_{U_i}$ satisfying
\begin{enumerate}
\item
for all $i\in I$,
$\Phi_i$ is $C_0(U_i)$-linear and $G$-equivariant, and
\item
for all $i,j\in I$, $b\in B$, and $x\in U_{ij}$,
$\Phi_j(b)(x)=\alpha^x_{\gamma_{ij}(x)}\big(\Phi_i(g)(x)\big)$.
\end{enumerate}
\end{prop}
\begin{proof}
Since $\gamma$ is a representative for $q:Z\to X$ in
$Z^1(X,\mathcal G)$, there exist local sections $\varphi_i:U_i\to q^{-1}(U_i)$
such that $\varphi_i(x)=\gamma_{ij}(x)\varphi_j(x)$ for all $x\in U_{ij}$.
It follows then from Remark~\ref{rem-trivial} that the isomorphisms
$\Phi_i:(\fc Z*A)_{U_i}\to A_{U_i}$ of Lemma~\ref{lem-trivial} corresponding
to the local sections $\varphi_i$ satisfy conditions \partref1~and \partref2.

Suppose now that $(B,G,\beta)$ is an arbitrary $C_0(X)$-system
and let $\Phi_i:B_{U_i}\to A_{U_i}$ be isomorphisms
satisfying
\partref1~and
\partref2. For each $z\in q^{-1}(U_i)$ define $s_i(z)\in G$ by the
equation
$z=s_i(z)\varphi_i\(q(z)\)$.  It follows from Lemma~\ref{lem-trivial}
that
\[
\Psi_i(b)(z)=
\alpha^{q(z)}_{s_i(z)^{-1}}\big(\Phi_i(b)(q(z))\big)
\]
defines a $C_0(U_i)$-linear and $G$-equivariant isomorphism $
\Psi_i:B_{U_i}\to (\fc Z*A)_{U_i}$ for all
$i\in I$.
Moreover, if $q(z)\in U_{ij}$, then
\[
z=s_i(z)\varphi_i\(q(z)\)=s_j(z)\gamma_{ij}\(q(z)\)\varphi_j\(q(z)\)
\]
which implies that $s_j(z)=s_i(z)\gamma_{ij}\(q(z)\)$ for all
$z\in Z$. It follows that
\begin{align*}
\Psi_i(b)(z)&=\alpha^{q(z)}_{s_i(z)^{-1}}\(\Phi_i(b)(q(z))\)
=\alpha^{q(z)}_{s_i(z)^{-1}}\(\alpha^{q(z)}_{\gamma_{ij}(q(z))^{-1}}
\(\Phi_j(b)(q(z))\)\)\\
&=\alpha^{q(z)}_{(s_i(z)\gamma_{ij}(q(z)))^{-1}}\(\Phi_j(b)(q(z))\)
=\alpha^{q(z)}_{s_j(z)^{-1}}\(\Phi_j(b)(q(z))\)\\
&=\Psi_j(b)(z),
\end{align*}
for all $z\in q^{-1}(U_{ij})$, $b\in B_{U_{ij}}$.
Thus, if we define $\Psi:B\to \fc Z*A$ by the formula
$\Psi(b)\(q(z)\)=\Psi_i(b)\(q(z)\)$, whenever
$z\in q^{-1}(U_i)$, it follows from Lemma~\ref{lem-useful}
that $\Psi$ is an isomorphism between the
$C_0(X)$-systems $(B,G,\beta)$ and $(\fc Z*A,G,\fc Z*\alpha)$.
\end{proof}

We close this section with the following useful result
which we will need later.
It will be helpful to keep in mind that if $A$ is a \coxalg, and $B$
any \cs-algebra, then \lemref{lem-maxtensor} implies that $B\mtensor A$
is a \coxalg{} with fibres $(B\mtensor A)_x\cong B\mtensor A_x$.

\begin{prop}\label{prop-tensor}
Let $q:Z\to X$ be a proper $G$-bundle with $G$ abelian.
Suppose further that $(A,G,\alpha)$ is a $C_0(X)$-system
and that  $B$ is a $C_0(X)$-algebra.
Then
$$\big(\fc Z*{(B\otimes_XA)},G,\fc Z*{(\id\otimes_X\alpha)}\big)\quad
\text{and}\quad \big(B\otimes_X(\fc Z*A),G,\id\otimes_X(\fc Z*\alpha)\big)$$
are isomorphic $C_0(X)$-systems.
\end{prop}
\begin{proof} We first show that for any system $(A,G,\alpha)$
and any $\cs$-algebra $B$ the induced
algebra $Z\times_G (B\mtensor A)$ is  equivariantly isomorphic to
$B\mtensor(\indgza)$. For this
let $(i_B,i_A)$ denote the natural embeddings
of $B$ and $A$ in $M(B\mtensor A)$, and define
homomorphisms
$\Phi_B:B\to M(Z\times_G (B\mtensor A))$ and
$\Phi_{\indgza}:\indgza\to M(Z\times_G (B\mtensor A))$ by
\[
(\Phi_B(b)\cdot H)(z)= i_B(b)(H(z))\quad\text{and}\quad
(\Phi_{\indgza}(F)\cdot H)(z)=i_A(F(z))H(z),
\]
where $b\in B$, $F\in \indgza$, and  $H\in Z\times_G (B\otimes_XA)$.
It is then straightforward to check that $\Phi_B$ and $\Phi_{\indgza}$
are commuting nondegenerate $*$-homomorphisms such that
$\set{\Phi_B(b)\Phi_{\indgza}(F): \text{$b\in B$ and $ F\in \indgza$}}$
generates a dense subalgebra of $Z\times_G (B\mtensor A)$.
Thus, by the universal property of the maximal tensor product we
obtain a surjective $*$-homomorphism
$\Phi_B\otimes\Phi_{\indgza}:B\mtensor(\indgza)\to
Z\times_G (B\mtensor A)$ which is clearly $G$-equivariant
and $C_0(X)$-linear with respect to the $C_0(X)$-structures
of the induced algebras. To see that $\Phi_B\otimes\Phi_{\indgza}$ is
an isomorphism, it suffices to see that the induced maps
\[
(\Phi_B\otimes\Phi_{\indgza})^x: (B\mtensor(\indgza))_x\to
(Z\times_G (B\mtensor A))_x
\]
are isomorphisms for all $x$.
To see this, note that both fibres
$(B\mtensor (\indgza))_x$ and $(Z\times_G (B\mtensor A))_x$
are isomorphic to $B\mtensor A$ (this follows from
Lemma~\ref{lem-maxtensor} in the first case).
If we do these identifications, then for $b\in B$ and $a\in A$
we can compute $(\Phi_B\otimes\Phi_{\indgza})^x(a\otimes b)$ as follows:
Choose $F\in \indgza$ and $z\in Z$ with $q(z)=x$ and $F(z)=a$.
Then
\[
(\Phi_B\otimes\Phi_{\indgza})^x(a\otimes b)=
\Phi_B\otimes\Phi_{\indgza}(b\otimes F)(z)= b\otimes F(z)=b\otimes a.
\]
Thus $(\Phi_B\otimes\Phi_{\indgza})^x$ is the identity on
$A\mtensor B$ for each $x\in X$.

Now let $(A,G,\alpha)$ and $B$ be as
in the proposition. Then
$B\mtensor (\indgza)$ and $Z\times_G (B\mtensor A)$ are
$C_0(X\times X\times X)$-algebras
(Lemma~\ref{lem-maxtensor}), and it follows directly from the
definition that
$\Phi_B\otimes\Phi_{\indgza}$ is $C_0(X\times X\times X)$-linear.
If $\Delta^3(X)=\set{(x,x,x):x\in X}$ denotes the diagonal in
$X\times X\times X$, then $B\otimes_X\fc Z*A$ is the restriction
of $B\mtensor (\indgza)$ to $\Delta^3(X)$ and
$\fc Z*(B\otimes_XA)$ is the restriction of $Z\times_G (B\mtensor A)$
to $\Delta^3(X)$. Thus the result follows from the
$C_0(X\times X\times X)$-linearity of $\Phi_B\otimes\Phi_{\indgza}$.
\end{proof}

%% file: inner.tex
%

\section{Crossed products by inner actions}\label{first-step}

In this section, we start with the investigation of crossed
products by inner actions on $\CR(X)$-algebras. If $\alpha:G\to
\Inn(A)$ is an inner action of a second countable group $G$
on a $\CR(X)$-algebra $A$, then it follows
from the discussion following \cite[Remark 2.9]{ew2}
that there exists a cocycle $u\in Z^2(G, C(X,\TT))$
and a Borel map
$v:G\to
\UM(A)$ satisfying
\begin{equation}\label{eq-cocycle}
 \alpha_s=\Ad v_s\quad\text{and}\quad v_sv_t=u(s,t)v_{st}
\quad\text{for all $s,t\in G$.}
\end{equation}
Then
\cite[Corollary 0.12]{rr} implies that the cohomology class $[u]\in
H^2(G,C(X,\TT))$ is a complete invariant for the
exterior equivalence class of $\alpha$.
(Throughout, $C(X,\T)$ is the trivial Polish $G$-module
$C(X,\TT)$ equipped with the topology of uniform convergence
on compact sets.)

To describe the crossed product $A\rtimes_{\alpha}G$ in terms of the
cocycle~$u$ it is necessary to work with (Busby-Smith) twisted actions
and twisted crossed products.  A {\em twisted action\/} $(\alpha,u)$
of a second countable locally compact group $G$ on a separable
\cs-algebra $A$ consists of a strongly measurable map $\alpha:G\to
\Aut(A)$ and a strictly measurable map $u:G\times G\to \UM(A)$
satisfying
\begin{enumerate}
\item
$\alpha_e=\id$ and $u(e,s)=u(s,e)=1$ for all $s\in G$;
\item
$\alpha_s(\alpha_t(a))=u(s,t)\alpha_{st}(a)u(s,t)^*$
for all $s,t\in G$;
\item
$\alpha_r(u(s,t))u(r,st)=u(r,s)u(rs,t)$
for all $s,t,r\in G$.
\end{enumerate}
The quadruple $(A,G,\alpha,u)$ is called a
(Busby-Smith or Leptin) {\em twisted
$\cs
$-dynamical system}. If $A$ is a $C_0(X)$-algebra, and $\alpha_s$ is
$C_0(X)$-linear for all $s\in G$, then we will call
$(A,G,\alpha,u)$ a {\em twisted $C_0(X)$-system}.

A {\em covariant homomorphism\/} of $(A,G,\alpha,u)$
into the multiplier algebra of a separable $\cs $-algebra $B$
is a pair
$(\Phi,v)$, where $\Phi:A\to M(B)$ is a nondegenerate
homomorphism and $v:G\to \UM(B)$ is strictly measurable
such that $v_e=1$, and such that
\[
\Phi(\alpha_s(a))=v_s\Phi(a)v_s^*,
\quad\text{and}\quad v_sv_t=\Phi(u(s,t))v_{st}
\quad\text{for all $s,t\in G$, $a\in A$.}
\]
 A
{\em covariant representation\/} of $(A,G,\alpha,u)$
on a separable Hilbert space $\H$ is a covariant
homomorphism $(\pi,U)$ of $(A,G,\alpha,u)$
into $\mathcal B(\H)$
(viewed as the multiplier algebra of $\K(\H)$).
If $(\Phi,v)$ is a covariant homomorphism, then
the integrated form $\Phi\rtimes v:L^1(G,A)\to M(B)$
is defined by
\[\Phi\rtimes v(f)=\int_G\Phi(f(s))v_s\,ds.\]
The following is a slight reformulation of
Packer and Raeburn's definition of a crossed product
for twisted systems.

\begin{definition}[cf., {\cite[Theorem
1.2]{para3}}]
\label{def-twistcross}
Let $(A,G,\alpha,u)$ be a
twisted system. A {\em crossed product} for $(A,G,\alpha,u)$
consists of a triple
$(B,i_A,i_G)$ satisfying
\begin{enumerate}
\item
$(i_A,i_G)$ is a covariant homomorphism of
$(A,G,\alpha,u)$ into $M(B)$;
\item
$i_A\rtimes i_G(L^1(G,A))$ is a dense subalgebra of $B$;
\item
if $(\Phi,v)$ is any covariant homomorphism
of $(A,G,\alpha,u)$ into $M(C)$, for some separable $\cs $-algebra $C$,
then there exists a nondegenerate $*$-homomorphism
$\Phi\rtimes v:B\to M(C)$ such that
$(\Phi\rtimes v)\circ i_A=\Phi$ and $(\Phi\rtimes v)\circ
i_G=v$.
\end{enumerate}
\end{definition}

If $(B,i_A,i_G)$ and $(C,j_A,j_G)$ are two different crossed products
of $(A,G,\alpha,u)$, then $j_A\rtimes j_G:B\to C$ is an isomorphism
with inverse $i_A\rtimes i_G:C\to B$.
Thus, the
crossed product is unique up to isomorphism, and we will usually
suppress the maps $i_A$ and $i_G$ and
denote it by
$A\rtimes_{\alpha,u}G$.

If $(A,G,\alpha,u)$ is a twisted $C_0(X)$-system,
then $A\rtimes_{\alpha,u}G$ is a $C_0(X)$-algebra,
where the action of $C_0(X)$ on $A\rtimes_{\alpha,u}G$
is defined via the composition of maps
\[
C_0(X)\arrow{e,t}{\phi}
\ZM(A)\arrow{e,t}{i_A}
M(A\rtimes_{\alpha,u}G).
\]
Moreover, the fibres $(A\rtimes_{\alpha,u}G)_x$ are
isomorphic  to $A_x\rtimes_{\alpha^x,u^x}G$,
where $(\alpha^x,u^x)$ is the twisted action induced
from $(\alpha,u)$ on the $G$-invariant quotient $A_x$
(see \cite{may1}
for more details).

For $\cox$-systems, the next observation follows immediately from
the equality $(\Phi\rtimes v)(i_A\circ \phi(f))=
\Phi(\phi(f))$ for all $f\in C_0(X)$.

\begin{lem}\label{lem-CXcov}
Let $(A,G,\alpha,u)$ be a twisted $C_0(X)$-system, and let $D$ be a
$C_0(X)$-algebra. If $(\Phi,v)$ is a covariant homomorphism
of $(A,G,\alpha,u)$ into $M(D)$, then $\Phi\rtimes v$ is
$C_0(X)$-linear if and only if $\Phi$ is $C_0(X)$-linear.
\end{lem}

If $(A,G,\alpha,u)$
is a twisted system, then there exists a
``dual'' action of the Pontryagin dual $\hgab$ of the
\emph{abelianization}
$\gab =G/\overline{[G,G]}$ of $G$
on $A\rtimes_{\alpha,u}G$.  This action is defined by
\begin{equation}
\label{eq-*dag}
\dual(\alpha,u)_{\chi}=i_A\rtimes (\bar\chi\cdot i_G),
\end{equation}
where $\bar\chi\cdot i_G(s)=\overline{\chi(s)}i_G(s)$
(here we view $\chi\in\hgab$  as a function on $G$),
and $i_A\rtimes (\bar\chi\cdot i_G)$
denotes the integrated form of the covariant homomorphism
$(i_A,\bar\chi\cdot i_G)$.
If $f\in L^1(G,A)$, viewed as a subspace of
$A\rtimes_{\alpha,u}G$, we have
$\dual(\alpha,u)_{\chi}(f)=\bar\chi f$, the pointwise
product of
$\bar\chi$ with $f$. If $(A,G,\alpha,u)$ is a twisted
$C_0(X)$-system, then $(A\rtimes_{\alpha,u}G, \hgab ,
\dual(\alpha,u) )$ is also a
$C_0(X)$-system since $i_A$ is $C_0(X)$-linear.
(Notice that
ordinary (separable)
$\cs $-dynamical systems and their crossed products are
recovered as the
special case where $u\equiv 1$, in which case we simply write
$(A,G,\alpha)$ for the system and $A\rtimes_{\alpha}G$
for the crossed product.)

The following proposition is crucial as it will allow us to
untangle certain diagonal twisted actions.
It is the $C_0(X)$-analogue
for twisted actions of the well known isomorphism of
$(A\mtensor B)\rtimes_{\alpha\otimes\id_B }G$ with
$(A\rtimes_{\alpha}G)\mtensor B$.

\begin{prop}\label{prop-twisttensor}
Let $(A,G,\alpha,u)$ be a twisted $C_0(X)$-system,
and let $B$ be a separable $C_0(X)$-algebra. Let $(\id,1)$ denote
the trivial $G$-action on $B$, and let $(\alpha\xtensor \id, u\xtensor 1)$
denote the diagonal twisted action of $G$ on $A\otimes_X B$.
If $i_A$ and $ i_G$ denote the canonical maps from
$A$ and $G$ into $M(A\rtimes_{\alpha,u}G)$,
then $(i_A\xtensor \id, i_G\xtensor 1)$ is a covariant
homomorphism of
$(A\xtensor B,G,\alpha\xtensor \id, u\xtensor  1)$
into $M((A\rtimes_{\alpha,u}G)\xtensor B)$. In particular, the
integrated form $(i_A\xtensor \id)\rtimes(i_G\xtensor 1)$ is a
$C_0(X)$-linear and $\hgab $-equivariant isomorphism
of $(A\xtensor B)\rtimes_{\alpha\xtensor \id, u\xtensor 1}G$ onto
$(A\rtimes_{\alpha,u}G)\xtensor B$.
\end{prop}

\begin{proof}
Put $C:=(A\rtimes_{\alpha,u}G)\xtensor B$.
To show that $(C, i_A\xtensor \id, i_G\xtensor 1)$ is a crossed
product for $(A\xtensor B, G, \alpha\xtensor \id, u\xtensor  1)$,
we have to verify conditions \partref1, \partref2, and \partref3
of Definition~\ref{def-twistcross}. Since $(i_A,i_G)$ is a
covariant homomorphism of
$(A,G,\alpha,u)$, it follows that
$(i_A\xtensor \id, i_G\xtensor  1)$ is a covariant homomorphism
of $(A\xtensor B,G,\alpha\xtensor  \id, u\xtensor  1)$.
This proves~\partref1.

For \partref2, let
$\Phi:L^1(G,A)\odot B\to L^1(G, A\xtensor B)$  be defined by
$\Phi(f\otimes b)(s) = f(s)\xtensor  b.$
Then
\begin{multline*}
(i_A\xtensor \id)\rtimes(i_G\xtensor 1)(\Phi(f\otimes b))
=
\int_G (i_A(f(s))\xtensor  b)(i_G(s)\xtensor 1)\,ds\\
=\left(\int_Gi_A(f(s))i_G(s)\,ds\right)\xtensor b
=\big(i_A\rtimes i_G(f)\big)\xtensor b.
\end{multline*}
Now \partref2 follows from the fact that
$i_A\rtimes i_G(L^1(G,A))$
is dense in $A\rtimes_{\alpha,u}G$.

For \partref3, suppose that $(\Psi, v)$ is a covariant
homomorphism of $(A\xtensor B,G,\alpha\xtensor  \id,
u\xtensor 1)$ into $M(D)$ for some separable $\cs $-algebra $D$.
By Remark~\ref{rem-universal}, we have
$\Psi=\Psi_A\xtensor \Psi_B$ such that
$\Psi_A(a\cdot f)\Psi_B(b)=\Psi_A(a)\Psi_B(f\cdot b)$ for
$f\in C_0(X)$, $a\in A$, and $b\in B$.
It is straightforward to check that
$(\Psi_A,v)$ is a covariant homomorphism of $(A,G,\alpha,u)$ into
$M(D)$ which commutes with $\Psi_B$. Thus we obtain a homomorphism
$(\Psi_A\rtimes v)\otimes
\Psi_B:(A\rtimes_{\alpha,u}G)\mtensor B\to M(D)$. For
$g\in L^1(G,A)$,
$b\in B$, and $f\in C_0(X)$ we have
\begin{align*}
\Psi_A\rtimes v (f\cdot g)\Psi_B(b)
&=\int_G\Psi_A(g(s)\cdot f) v(s)\,ds\Psi_B(b)
=\int_G\Psi_A(g(s)\Psi( f)  v(s)\, ds\Psi_B(b)\\
&=\int_G\Psi_A(g(s))v(s)\,ds\Psi_B(f\cdot b)
=\Psi_A\rtimes v(g)\Psi_B(f\cdot b).
\end{align*}
Since $L^1(G,A)$ is dense in $A\rtimes_{\alpha,u}G$ this extends to all
$g\in A\rtimes_{\alpha,u}G$.
By Remark~\ref{rem-universal}, there is a nondegenerate
homomorphism $(\Psi_A\rtimes v)\xtensor
\Psi_B: C\to M(D)$  satisfying $(\Psi_A\rtimes v)\xtensor
\Psi_B(g\xtensor  b)=\Psi_A\rtimes v(g)\Psi_B(b)$ for all
elementary tensors $g\xtensor b$.

Finally, we compute
\begin{align*}
(\Psi_A\rtimes v)\xtensor  \Psi_B(i_A\xtensor \id(a\xtensor b))
&=(\Psi_A\rtimes v)(i_A(a))\Psi_B(b)\\
&=\Psi_A(a)\Psi_B(b)=\Psi(a\xtensor b),
\end{align*}
which proves that $(\Psi\rtimes v)\circ(i_A\xtensor \id)=\Psi$.
Similarly,
\[
(\Psi_A\rtimes v)\xtensor  \Psi_B(i_G(s)\xtensor 1) =
(\Psi_A\rtimes v)(i_G(s))=v_s,
\]
which proves that $(\Psi\rtimes v)\circ (i_G\xtensor 1)=v$, as
required.

At this point, we have proved that $(C, i_A\xtensor \id,i_G\xtensor 1)$
is a crossed  product for $(A\xtensor B, G,\alpha\xtensor \id, u\xtensor
1)$. In particular,
\[(i_A\xtensor \id)\rtimes(i_G\xtensor 1):
(A\rtimes_{C_0(X)}B)\rtimes_{\alpha\xtensor \id, u\xtensor 1}G\to
(A\rtimes_{\alpha,u}G)\xtensor B\] is an isomorphism,
which is $C_0(X)$-linear since $i_A\xtensor\id$ is $C_0(X)$-linear.
If $\chi\in \hgab $, then
$\bar\chi\cdot (i_G\xtensor 1)=(\bar\chi\cdot i_G)\xtensor 1$,
which implies that $(i_A\xtensor \id)\rtimes(i_G\xtensor 1)$ also preserves
the dual action of $\hgab $.
\end{proof}

\begin{definition}
[{\cite[Definition~3.1]{para1}}]
\label{def-exteq}
Two twisted actions $(\alpha, u)$ and $(\beta,v)$
of $G$ on $A$ are called
{\em exterior equivalent} if there exists a strictly Borel map
$w:G\to
\UM(A)$ satisfying
\[
\beta_s=\Ad w_s\circ \alpha_s\quad\text{and} \quad
v(s,t)=w_s\alpha_s(w_t)u(s,t)w_{st}^*\quad\text{for all $s,t\in G$.}
\]
\end{definition}

Note that if $(A,G,\alpha,u)$ is a twisted $C_0(X)$-system,
and if $(\beta,v)$ is exterior equivalent to $(\alpha,u)$,
then
\[
\beta_s(f\cdot a)=w_s\alpha_s(f\cdot a)w_s^*=
f\cdot (w_s\alpha_s(a)w_s^*)=f\cdot \beta_s(a),
\]
so that each $\beta_s$ is $C_0(X)$-linear.
Further, if  $w_s^x$ denotes the image of $w_s$ in
$M(A_x)$,  then
$w^x$ implements an exterior equivalence between
$(\alpha^x,u^x)$ and $(\beta^x, v^x)$.

\begin{lem}
[{\cite[Lemma 3.3]{para1}}]
\label{lem-exisom}
Suppose that $w:G\to \UM(A)$ implements an exterior equivalence
between the twisted actions $(\alpha,u)$ and $(\beta,v)$ of $G$ on
$A$. Let $j_A:A\to M(A\rtimes_{\beta,v}G)$ and $j_G:G\to
\UM(A\rtimes_{\beta,v}G)$ denote the canonical maps and let
$\mu_G:G\to \UM(A\rtimes_{\beta,v}G)$ be defined by
$\mu_G(s)=j_A(w_s^*)j_G(s)$. Then
$(j_A,\mu_G)$ is a covariant homomorphism of
$(A,G,\alpha, u)$, and $j_A\rtimes \mu_G$
is an isomorphism between $A\rtimes_{\alpha,u}G$ and
$A\rtimes_{\beta,v}G$.
\end{lem}

\begin{remark}\label{rem-exisom}
If $(A,G,\alpha,u)$ is a twisted $C_0(X)$-system and
$(\beta,v)$ is exterior equivalent to $(\alpha,u)$ via
$w:G\to \UM(A)$, then the isomorphism $j_A\rtimes \mu_G$ above is
necessarily
\cox-linear and therefore
implements an isomorphism between the $C_0(X)$-systems
$(A\rtimes_{\alpha,u}G,\hgab ,\dual(\alpha,u) )$
and
$(A\rtimes_{\beta,v}G,\hgab ,\dual(\beta,v))$.
To see this, note that $j_A$ is $C_0(X)$-linear,
which implies $C_0(X)$-linearity of $j_A\rtimes \mu_G$ (\lemref{lem-CXcov}).
Further, if $\chi\in \hgab $, then it follows from
the definition of the dual actions that
\begin{equation}\label{eq-new}
\begin{split}
(j_A\rtimes \mu_G)\circ \dual(\alpha,u)_{\chi}&=
j_A\rtimes (\bar\chi\cdot \mu_G)
=j_A\rtimes \big(\bar\chi\cdot ((j_A\circ w)\cdot j_G)\big)\\
&=j_A\rtimes \big((j_A\circ w)\cdot(\bar\chi\cdot j_G)\big)
=\dual(\beta,v)_{\chi}\circ (j_A\rtimes \mu_G).
\end{split}
\end{equation}
\end{remark}

Twisted crossed products with $A=C_0(X)$ abelian and $\alpha$
trivial play a central r\^ole in our analysis.
Such a crossed product is called a twisted
transformation group \cs-algebra and is denoted by $\cs(G,X,u)$;
a
nice survey article is
\cite{judy-exp}.  For our purposes, we need only remark that the
condition on the twist
$u$ implies that
$u$ is a cocycle in the Moore cohomology group $Z^2(G,C(X,\TT))$
for the trivial action of $G$ on $C(X,\TT)$. Moreover, if $u, v\in
Z^2(G,C(X,\TT))$, then $(\id,u)$ is exterior equivalent to $(\id,v)$
if and only if $[u]=[v]$ in $H^2(G,C(X,\TT))$. In
what follows we shall  write
$\hatu $ (instead of $\dual(\id,u)$ as in \eqref{eq-*dag})
for the dual action of $\hgab $ on $\cs (G,X,u)$. It follows
from Lemma~\ref{lem-exisom} that $(\cs(G,X,u),\hgab,\hat{u})$
only depends on the cohomology class of~$u$.

If $A$ is a \coxalg, then since the action of $\cox$ on $A$ is
nondegenerate, we can extend the action of $\cox$ to one of $C^b(X)$.
Therefore, we can make the following definition.

\begin{definition}\label{def-u-homomorphism}
Suppose that $A$ is a separable $C_0(X)$-algebra, and that
$u\in Z^2(G,C(X,\TT))$. A {\em $u$-homomorphism\/} is a strictly
measurable map
$v:G\to \UM(A)$  satisfying
\[v_e=1\quad\text{and}\quad v_sv_t=u(s,t)\cdot v_{st}
\quad\text{for all $s,t\in G$.}
\]
If $\alpha:G\to \Aut(A)$ is an action, then we say that $\alpha$
is {\em implemented\/} by the $u$-homomorphism $v:G\to \UM(A)$
if $\alpha_s=\Ad v_s$ for all $s\in G$.
\end{definition}

Notice that if $v$ is a $u$-homomorphism and if $\phi:C_0(X)\to \ZM(A)$
is the homomorphism determined by the $C_0(X)$-action on $A$, then
$(\phi,v)$ is a covariant homomorphism of $(C_0(X),G,\id,u)$ into
$M(A)$. We will write
$\bar u$ for the inverse of the cocycle $u\in Z^2(G,C(X,\TT))$.
If $\alpha:G\to \Inn(A)$ is implemented by the $u$-homomorphism
$v:G\to \UM(A)$,  then it follows from
$\alpha_s=\Ad v_s$ and $1=v_sv_t\bar{u}(s,t)v_{st}^*$ for $s,t\in G$,
that $v$ implements an exterior equivalence between the twisted actions
$(\alpha,1)$ and $(\id,\bar u)$. We use this observation
for the proof of

\begin{prop}\label{prop-twistdecom}
Suppose that  $A$ is a \coxalg, that $u\in Z^2(G,C(X,\TT))$, and that
$\alpha:G\to\Aut(A)$ is implemented by a $u$-homomorphism
$v:G\to \UM(A)$. Then
$A\rtimes_{\alpha}G$ is isomorphic to $\cs (G,X,\bar{u})\xtensor A$.
In particular, if  $i_G:G\to \UM(\cs (G,X,\bar{u}))$
is the canonical map, then $(1\xtensor \id_A, i_G\xtensor v)$
is a covariant homomorphism of $(A,G,\alpha)$ into
$M(\cs (G,X,\bar u)\xtensor A)$ whose
integrated form is a $C_0(X)$-linear covariant isomorphism of
$(A\rtimes_{\alpha}G,\hgab, \widehat{\alpha})$ onto
$(\cs (G,X,\bar u)\xtensor A,\hgab,\widehat{\bar u} \xtensor~\id)$.
\end{prop}
\begin{proof} It follows from Proposition~\ref{prop-twisttensor} that
$(i_{C_0(X)}\xtensor  \id_A,  i_G\xtensor 1)$ is a
covariant homomorphism of
$(C_0(X)\xtensor A, G, \id\xtensor \id, \bar{u}\xtensor  1)$ into
$M(\cs (G,X, \bar u)\xtensor A)$ whose integrated form
is a $C_0(X)$-linear and $\hgab $-equivariant isomorphism
of $(C_0(X)\xtensor A)\rtimes_{\id\otimes\id,  \bar u\otimes 1}G$
onto $\cs (G,X,\bar u)\xtensor A$.

Let $\Phi:C_0(X)\xtensor A\to A$ be the isomorphism defined
on elementary tensors by $\Phi(f\xtensor a)=f\cdot a$.
Then $\Phi$ carries the trivial action $\id\otimes \id$
to the trivial action on $A$
and we have $\Phi(\bar u(s,t)\xtensor 1)a=\bar u(s,t)\cdot a$ for all
$a\in A$. Thus, regarding $\bar u$ as a map $\bar u:G\times G\to \UM(A)$
via the
$C_0(X)$-action on $A$, we see that $\Phi$ induces
a $C_0(X)$-linear and $\hgab $-equivariant isomorphism between
$(C_0(X)\xtensor A)\rtimes_{\id\xtensor \id, \bar{u}\xtensor  1}G$
and $A\rtimes_{\id,\bar u}G$.
Moreover, $\Phi$ carries the covariant homomorphism
$(i_{C_0(X)}\xtensor \id_A, i_G\xtensor 1)$ to the
covariant homomorphism $(1\xtensor \id_A, i_G\xtensor 1)$, which implies
that
$(1\xtensor \id_A)\rtimes(i_G\xtensor 1)$ is a $C_0(X)$-linear and
$\hgab $-equivariant homomorphism of
$A\rtimes_{\id,\bar u}G$ onto $\cs (G,X,\bar u)\xtensor A$.
Since $v$ is a
$u$-homomorphism, it follows that $v$ implements an exterior
equivalence  between the twisted actions $(\alpha,1)$ and
$(\id, \bar u)$. Thus the result follows from Lemma~\ref{lem-exisom}.
\end{proof}

If
$\alpha:G\to\Inn(A)$ is any inner action of $G$ on the $\CR(X)$-algebra
$A$, then as we observed at the beginning of this section,
there exists a
unique class
$[u]\in H^2(G,C(X,\TT))$ and a $u$-homomorphism $v:G\to \UM(A)$ which
implements $\alpha$. Thus we get

\begin{cor}\label{cor-inner}
Let $A$ be a $\CR(X)$-algebra and let $\alpha:G\to\Inn(A)$
be an inner action of $G$ on $A$. Let $u\in Z^2(G,C(X,\TT))$
be associated to $\alpha$ as above. Then
$(A\rtimes_{\alpha}G,\hgab,\widehat{\alpha})$ is
$C_0(X)$-isomorphic to $(\cs(G,X,\bar u)\otimes_X A,
\hgab,\widehat{\bar u}\otimes_X\id)$.
\end{cor}

\begin{remark}\label{rem-inner}
It is shown in
\cite[Proposition 3.1]{horr} that for every $u\in Z^2(G,C(X,\TT))$,
there exists a $u$-homomorphism $v$ of $G$ into $\UM\(C_0(X,\K)\)$.
From this it follows that for any stable $C_0(X)$-algebra $A$,
there exists a $u$-homomorphism $w:G\to \UM(A)$: simply identify
$A\cong A\otimes\K$ with
$A\otimes_XC_0(X,\K)$ and define
$w=1\otimes_Xv$. Thus, if $A\in \CR(X)$ is stable, then
there exists a natural one-to-one correspondence between the
exterior equivalence classes of inner $G$-actions on $A$ and
$H^2(G,C(X,\TT))$, and by the above results the crossed products
can be described in terms of the central twisted
transformation group algebras $\cs(G,X,\bar u)$.
\end{remark}

%% file: loc-fibre.tex
%
%
\section{$\hgab$-fibre products and locally unitary actions on
$C_0(X,\K)$}\label{loc-fibre}

The exterior equivalence classes of locally unitary actions
on $A$  are classified by the isomorphism classes
of principal $\hgab$-bundles, or equivalently,
 by classes in
$H^1(X,\shgab)$ as described in
\secref{sec-bund-op.6}.
This correspondence was originally worked out (for
abelian groups acting on continuous-trace \cs-algebras) by Phillips and
Raeburn \cite{pr2}; the details of the
extension to arbitrary groups acting on
\CR-algebras can be found in \cite[\S3]{ew2}.
Recall that if $A\in \CR(X)$, then an action $\alpha:G\to \Aut(A)$
is called  {\em locally unitary\/} if each point in $X$
has an open neighborhood $W$ such that
the restriction $\alpha^W$ of $\alpha$ to the ideal $A_W$ of $A$
is unitary; that is, $\alpha^W=\Ad w$ for some
strictly continuous homomorphism $w:G\to \UM(A_W)$.
The class corresponding to a
locally unitary action is determined as follows.
Choose any locally finite open cover $(W_i)_{i\in I}$ of $X$ such that
each restriction $\alpha^i:=\alpha^{W_i}$ is unitary.
For each $i\in I$, set $A_i:=A_{W_i}$ and let $w^i:G\to \UM(A_i)$
  be a strictly
continuous map such that $\alpha^i=\Ad w^i$. If $w^i(s,x)$ denotes the
element of $M(A_x)$ induced by $w^i$, then
there exist continuous functions
$\gamma_{ij}:W_{ij}\to \hgab$ satisfying
\begin{enumerate}
\item
$w^i(s,x)=\gamma_{ij}(x)(s)w^j(s,x)$ for all $x\in W_{ij}:=W_i\cap W_j$,
\item
$\gamma_{ij}(x)\gamma_{jk}(x)=\gamma_{ik}(x)$ for all $x\in
W_{ijk}:=W_i\cap W_j\cap W_k$.
\end{enumerate}
The last property implies that $(\gamma_{ij})_{i,j\in I}$
is a cocycle in $Z^1(X, \shgab)$, and by \cite[Proposition 3.3]{ew2}
the class
$\zeta(\alpha)$ of this cocycle in $H^1(X,\shgab)$ is a complete invariant
for the  exterior equivalence class of $\alpha$. If $A$ is stable,
then all classes in $H^1(X,\shgab)$ appear this way.

\begin{remark}\label{rem-stable}
Since $\coxk$ is stable, given  $\zeta\in H^1(X,\shgab)$, there
is a locally unitary action
$\delta:G\to \Aut(C_0(X,\K))$ with $\zeta(\delta)=\zeta$. Thus, if
$\alpha:G\to \Aut(A)$ is a locally unitary action on a
$\CR(X)$-algebra $A$, then
there exists a locally unitary action $\delta:G\to \Aut(C_0(X,\K))$
with $\zeta(\delta)=\zeta(\alpha)$.  Then \cite[Lemma 3.5]{ew2}
implies that, after identifying $A\tensor\K$ with $A\otimes_XC_0(X,\K)$,
$\zeta(\alpha\otimes\id_{\K})=\zeta(\id_A\otimes_X\delta)$.
In particular, $\alpha\otimes\id_{\K}$ is exterior equivalent to
$\id_A\otimes_X\delta$.
\end{remark}

Now let  $(A,G,\alpha,u)$ be any twisted $C_0(X)$-system.
If $\delta:G\to \Aut(C_0(X,\K))$ is locally unitary,
then we want to describe
the crossed product with respect to the twisted diagonal action
$(\alpha\otimes_X\delta, u\otimes_X1)$ of $G$ on
$A\otimes_XC_0(X,\K)\cong A\otimes \K$ in terms of
$A\rtimes_{\alpha,u}G$ and a
principal $\hgab$-bundle $q:Z\to X$ corresponding to $\zeta(\delta)$.
In particular, using Remark~\ref{rem-stable}, we will obtain
a description of $A\rtimes_{\alpha}G$ for any locally unitary
action $\alpha$ in terms of the bundle corresponding to $\zeta(\alpha)$.
But our more general result will be crucial for the description
of crossed products by locally inner actions given in the next section.
In order to state our result we need the following definition.

\begin{definition}\label{def-locext}
Let $(\alpha,u)$ and $(\beta,v)$ be two twisted $C_0(X)$-linear
actions of $G$ on a $C_0(X)$-algebra $A$. Then we say that
$(\alpha,u)$  is {\em locally exterior equivalent} to $(\beta,v)$
if every point in $X$ has an open neighborhood $W$
such that $(\alpha^W, u^W)$ is exterior equivalent to $(\beta^W, v^W)$.
\end{definition}

\begin{remark}\label{rem-locuni}
If $A\in \CR(X)$, then it follows directly from the definitions that
$\alpha:G\to\Aut(A)$ is locally unitary if and only if $\alpha$ is
locally exterior equivalent to the trivial action $\id_A$.
\end{remark}

\begin{thm}\label{thm-mostimportant}
Let $(\alpha,u)$ and $(\beta,v)$ be two locally exterior equivalent
twisted $C_0(X)$-linear actions of $G$ on a $C_0(X)$-algebra $A$.
Suppose further that $\delta:G\to\Aut(C_0(X,\K))$ is locally unitary
and that $q:Z\to X$ is a principal $\hgab$-bundle
corresponding to $\zeta(\delta)\in
H^1(X,\shgab)$. If $(\beta\otimes\id,v\otimes 1)$
is exterior equivalent to
$(\alpha\otimes_X\delta, u\otimes_X1)$ as actions on $A\otimes\K\cong
A\otimes_XC_0(X,\K)$, then
$\(A\rtimes_{\beta,v}G, \hgab,\dual(\beta, v)\)$
is $C_0(X)$-isomorphic to the $G$-fibre product
$\(\fc Z*{(A\rtimes_{\alpha,u}G)},\hgab, \fc Z*{\dual(\alpha,u)}\)$.
\end{thm}

Before we start with the proof we need the following two lemmas.
The first one follows directly from the formulas for
exterior equivalences given in Definition~\ref{def-exteq}.

\begin{lem}\label{lem-ext}
Let $(\alpha,u)$ be a twisted action of $G$ on $A$ and let
$\ZUM(A)=\UM(A)\cap \ZM(A)$.
Then a strictly measurable map $\lambda:G\to \UM(A)$ implements
an exterior equivalence of $(\alpha,u)$ with itself if and only
if $\lambda_s\in \ZUM(A)$  and $\lambda_{st}=
\lambda_s\alpha_s(\lambda_t)$ for all $s,t\in G$.
\end{lem}

\begin{lem}\label{lem-centralisom}
Let $A$ be a $\cs$-algebra. Then the isomorphism
$\ZM(A)\to \ZM(A\otimes \K)$ given by $z\mapsto z\otimes 1$
induces a homeomorphism between $\ZUM(A)$ and $\ZUM(A\otimes \K)$
with respect to the strict topologies.
\end{lem}
\begin{proof} First note that it follows from the
Dauns-Hofmann Theorem that the map $z\mapsto z\otimes 1$ is indeed an
isomorphism between $\ZM(A)$ and $\ZM(A\otimes \K)$.
If $(z_i)_{i\in I}$ converges strictly to $z$ in
$\ZUM(A)$, $z_ia\otimes c\to za\otimes c$ for every elementary
tensor $a\otimes c\in A\otimes \K$. Since
$\ZUM(A\otimes \K)$ is a bounded subset of $\M(A\otimes \K)$
this is enough to prove that $z_i\otimes 1\to z\otimes 1$ in
$\ZUM(A\otimes\K)$. Conversely, if  $z_i\otimes 1\to z\otimes 1$
strictly in $\ZUM(A)$, then choose any $c\in \K$ with $\|c\|=1$
to deduce that $\|z_ia-za\|=\|z_ia\otimes c-za\otimes c\|\to 0$
for all $a\in A$.
\end{proof}

\begin{remark}\label{rem-crossed}
If $(\alpha,u)$ is a twisted $C_0(X)$-linear
action of $G$ on $A$, and if $(j_A,j_G)$ denote the natural embeddings of $A$
and $G$ in $\M(A\rtimes_{\alpha,u}G)$, then for any nonempty open subset $W$ of
$X$ one can deduce from the definition of a crossed product that
$((A\rtimes_{\alpha,u}G)_W, j_A^W, j_G^W)$ is a crossed product for
$A_W\rtimes_{\alpha^W,u^W}G$, where
$j_A^W$ denotes the restriction of $j_A$ to the ideal $A_W$ of $A$
composed with the natural map
$\Phi:\M(A\rtimes_{\alpha,u}G)\to \M((A\rtimes_{\alpha,u}G)_W)$, and
$j_G^W=\Phi\circ j_G$. Similarly, for each
$x\in X$, a crossed product for $(A_x,G,\alpha^x,u^x)$ is given by
$((A\rtimes_{\alpha,u}G)_x, j_A^x,j_G^x)$, where
$j_A^x$ and $j_G^x$ are the compositions of $j_A$ and $j_G$
with the quotient map $\M(A\rtimes_{\alpha,u}G)\to
\M((A\rtimes_{\alpha,u}G)^x)$. Moreover,
for any $x\in W$ we have $(j_A^x\times j_G^x)\circ (j_A^W\times j_G^W)=
j_A^x\times j_G^x$.
\end{remark}

\begin{proof}[Proof of Theorem~\ref{thm-mostimportant}]
Since the
isomorphism class of
$\(\fc Z*{(A\rtimes_{\alpha,u}G)},\hgab, \fc Z*{\dual(\alpha,u)}\)$
only depends on the isomorphism class of $q:Z\to X$,
it follows from the
assumptions and the discussion at the beginning of this section
that we can  find a locally finite open cover
$(W_i)_{i\in I}$ of $X$ and a cocycle $(\gamma_{ij})_{i,j\in I}$ in
$Z^1(X,\shgab)$ (with respect to this cover) satisfying the
conditions:
\begin{enumerate}
\item
There exist strictly continuous maps $w^i:G\to \UM(C_0(W_i,\K))$
such that $\delta^i=\Ad w_i$ for all $i\in I$,
and $w^i(s,x)=\gamma_{ij}(x)(s)w^j(s,x)$ for all $x\in W_{ij}$.
\item
There exist continuous local sections $\varphi_i:W_i\to
q^{-1}(W_i)$ such that $\varphi_i(x)=\gamma_{ij}(x)\varphi_j(x)$
for all $x\in W_{ij}$.
\item
There exist strictly measurable maps $\kappa^i:G\to \UM(A_i)$,
$A_i:=A_{W_i}$, implementing exterior equivalences between
$(\beta^i,v^i)$ and $(\alpha^i, u^i)$; that is, $\beta^i=\Ad \kappa^i\circ
\alpha^i$ and
$v^i(s,t)=\kappa^i_s\alpha^i_s(\kappa^i_t)u^i(s,t)(\kappa^i_{st})^*$
for all $s,t\in G$.
\end{enumerate}
We want to use these data to define $C_0(W_i)$-linear and
$\hgab$-equivariant isomorphisms
$\Phi^i: (A\rtimes_{\beta,v}G)_{W_i}\to (A\rtimes_{\alpha,u}G)_{W_i}$
which satisfy
\begin{equation}\label{eqtrans}
\Phi^j(d)(x)=\dual(\alpha^x,u^x)_{\gamma_{ij}(x)}(\Phi^i(d))(x)
\quad\text{for all $d\in A\rtimes_{\beta,v}G$ and $x\in W_{ij}$.}
\end{equation}
If this can be done, then since  the bundle $q:Z\to
X$ has
transition functions $(\gamma_{ij})_{i,j\in I}$, the result follows from
Proposition~\ref{prop-principal} and the observation that the
action on $A_x\rtimes_{\alpha^x,u^x}G$ induced by $\dual(\alpha,u)$ is
$\dual(\alpha^x,u^x)$.

We claim that we may assume that the exterior equivalences
$\kappa_i:G\to \UM(A_i)$
satisfy the relation
\begin{equation}\label{eq-nu}
\kappa_i(s,x)=\gamma_{ij}(x)(s)\kappa_j(s,x)
\quad\text{for all $i,j\in I$, $x\in W_{ij}$ and $s\in G$.}
\end{equation}
To make notation easier we think of $(\alpha\otimes_X\delta, u\otimes_X1)$
as the family
of actions $(\alpha^x\otimes\delta^x, u^x\otimes 1)_{x\in X}$,
and we identify $A_i\tensor\K$ with $(A\xtensor C_0(X,\K))_{W_i}\cong
A_i\tensor_{W_i} C_0(W_i,\K)$.
Then (by committing a criminal abuse of notation) we denote
the restriction of $(\alpha\otimes_X\delta, u\otimes_X1)$ to
$A_i\otimes\K$
by $(\alpha^i\otimes\delta^i,u^i\otimes 1)$.
By assumption there exists a strictly measurable map
$\mu:G\to \UM(A\otimes \K)$
which implements an exterior equivalence between
$(\beta\otimes \id, v\otimes 1)$
and
$(\alpha\otimes_X\delta, u\otimes_X1)$, and we denote by $\mu^i$
the restriction of $\mu$ to $A_i\otimes\K$.
Since on each $W_i$ the map $w^i:G\to \UM(C_0(W_i,\K))$
implements an exterior equivalence between $\delta^i$ and the trivial
action on $C_0(W_i,\K)$, we can combine this with $\mu^i$
in order to obtain exterior equivalences $\sigma^i$ between
$(\beta^i\otimes\id, v^i\otimes 1)$ and $(\alpha^i\otimes \id,
u^i\otimes 1)$
given by
\begin{equation}\label{eq-sigma}
\sigma^i(s,x)=\mu^i(s,x)(1\otimes w^i(s,x)),\quad (s,x)\in G\times W_i.
\end{equation}
At this point we have two exterior equivalences between
$(\beta^i\otimes\id, v^i\otimes 1)$ and $(\alpha^i\otimes\id,
u^i\otimes 1)$, namely $\sigma^i$ and $\kappa^i\otimes 1$.
But then $s\mapsto \lambda^i(s)=(\kappa^i(s)\otimes 1)^*\sigma^i(s)$ is
an exterior equivalence for $(\alpha^i\otimes 1, u^i\otimes 1)$
with itself. Thus, Lemma~\ref{lem-ext} implies that
$\lambda^i$ takes values in $\ZUM(A_i\otimes \K)$ and satisfies
$\lambda^i(st)=\lambda^i(s)(\alpha^i\otimes\id)_s(\lambda^i(t))$.
It follows from Lemma~\ref{lem-centralisom} that there exists a
strongly measurable map $\tilde{\lambda}^i:G\to \ZUM(A_i)$ such
that $\lambda^i(s)=\tilde{\lambda^i}(s)\otimes 1$ for all $s\in G$, and
which satisfies
$\tilde{\lambda}^i(st)=\tilde{\lambda}^i(s)\alpha^i_s(\tilde{\lambda}^i(t))$
for all $s,t\in G$.
Therefore, $\tilde{\lambda}^i$ is an exterior equivalence
between $(\alpha^i,u^i)$ and itself.
Thus,
\[\tilde{\kappa}^i(s)=\kappa^i(s)\tilde{\lambda}^i(s)
\quad\text{for all $s\in G$}\]
defines a new exterior equivalence between $(\beta^i, v^i)$
and $(\alpha^i, u^i)$. Since $\mu^i(s,x)=\mu^j(s,x)$
if
$x\in W_{ij}$, we get
\begin{align*}
\tilde{\kappa}^i(s,x)\tilde{\kappa}^j(s,x)^*\otimes 1&
=\kappa^i(s,x)\tilde{\lambda}^i(s,x)\tilde{\lambda}^j(s,x)^*
\kappa^j(s,x)^*\otimes 1\\
&=\big(\kappa^i(s,x)\otimes
1\big)\lambda^i(s,x)\lambda^j(s,x)^*\big(\kappa^j(s,x)^*\otimes
1\big)\\
\intertext{which, by \eqref{eq-sigma}, is}
&=\sigma^i(s,x)\sigma^j(s,x)^*
=\mu^i(s,x)(1\otimes w^i(s,x))(1\otimes w^j(s,x)^*)\mu^j(s,x)^*\\
&=\mu^i(s,x)(1\otimes w^i(s,x)w^j(s,x)^*)\mu^j(s,x)^*\\
&=\gamma_{ij}(x)(s)\mu^i(s,x)\mu^j(s,x)^*
=\gamma_{ij}(x)(s)\otimes 1.
\end{align*}
Thus we see that the
$\tilde{\kappa}^i(x,s)=\gamma_{ij}(x)(s)\tilde{\kappa}^j(x,s)$ for all $i,j\in
I$ and $x\in W_{ij}$.
Thus we can replace $\kappa^i$ by $\tilde{\kappa}^i$ so that
\eqref{eq-nu} holds.  This proves the claim.

We now identify $(A\rtimes_{\beta,v}G)_{W_i}$ and
$(A\rtimes_{\alpha,u}G)_{W_i}$ with $A_i\rtimes_{\beta^i,v^i}G$ and
$A_i\rtimes_{\alpha^i, u^i}G$ as in Remark~\ref{rem-crossed},
and we let $j_A^i=j_A^{W_i}$ and $j_G^i=j_G^{W_i}$
denote the natural embeddings
of $A_i$ and $G$ in $\M(A\rtimes_{\alpha^i,u^i}G)$.
Then it follows from Lemma~\ref{lem-exisom} and Remark~\ref{rem-exisom}
that
\[
\Phi^i:=j_A^i\rtimes (j_A^i\circ (\kappa^i))\cdot
j_G^i:A_i\rtimes_{\beta^i,v^i}G
\to A_i\rtimes_{\alpha^i,u^i}G
\]
is indeed a $C_0(W_i)$-linear and $\hgab$-equivariant isomorphism for
all $i\in I$. Moreover, for all $x\in W_{ij}$ and $d\in
A\rtimes_{\beta,v}G$ we get (again using Remark~\ref{rem-crossed})
\begin{align*}
\Phi^j(d)(x)&=\big(j_A^x\rtimes (j_A^x\circ \kappa^j(\cdot,x))\cdot
j_G^x\big)(d(x))\\
&=\big(j_A^x\rtimes \gamma_{ij}(x)^{-1}\cdot (j_A^x\circ
\kappa^i(\cdot,x))\cdot j_G^x\big)(d(x))\\
&=\dual(\alpha^x,u^x)_{\gamma_{ij}(x)}\big(\big(j_A^x\rtimes
(j_A^x\circ \kappa^i(\cdot,x))\cdot j_G^x\big)(d(x))\big)\\
&=\dual(\alpha^x,u^x)_{\gamma_{ij}(x)}\big(\Phi^i(d)(x)\big).
\end{align*}
Thus the $\Phi^i$ satisfy equation (\ref{eqtrans}) which completes the proof.
\end{proof}

\begin{cor}\label{cor-twist}
Let $(\alpha,u)$ be a twisted $C_0(X)$-linear action of $G$ on $A$ and
let $\delta:G\to\Aut(C_0(X,\K))$ be locally unitary.
Let $q:Z\to X$ be a principal $\hgab$-bundle corresponding to
$\zeta(\delta)$. Then there is a
$C_0(X)$-linear and
$\hgab$-equivariant isomorphism between
$(A\otimes_XC_0(X,\K))\rtimes_{\alpha\otimes_X\delta, u\otimes_X1}G$
and $\fc Z *{(A\rtimes_{\alpha,u}G)}\otimes\K$.
\end{cor}
\begin{proof}
First recall that there exists a $C_0(X)$-linear and $\hgab$-equivariant
isomorphism between
$ (A\rtimes_{\alpha,u}G)\otimes\K$  and
$(A\otimes\K)\rtimes_{\alpha\otimes\id, u\otimes 1}G$,
and it follows from Proposition~\ref{prop-tensor} that this isomorphism
induces a $C_0(X)$-linear and $\hgab$-equivariant
isomorphism between
$\fc Z*{ (A\rtimes_{\alpha,u}G)}\otimes\K$  and
$\fc Z*{\big((A\otimes\K)\rtimes_{\alpha\otimes\id, u\otimes 1}G\big)}$
(with respect to the $\hgab$-actions
$\fc Z*{\dual(\alpha,u)}\otimes\id$
and $\fc Z*{(\specnp{(\alpha\otimes\id, u\otimes 1)})}$).
Since $(\alpha\otimes_X\delta, u\otimes_X1)$ is
locally exterior equivalent to $(\alpha\otimes\id, u\otimes 1)$
we can apply Theorem~\ref{thm-mostimportant} and the result follows.
\end{proof}

With \thmref{thm-mostimportant} in hand, we can give a general
description of crossed products by
locally unitary actions.  (Of course, if we were only interested in
locally unitary actions, we could have achieved this by more direct
methods.)

\begin{thm}\label{thm-loc-uni}
Let $A\in \CR(X)$ and let $\alpha:G\to \Aut(A)$
be locally unitary. Let $q:Z\to X$ be a principal
$\hgab$-bundle corresponding to $\zeta(\alpha)$, and let
$\mu:\hgab\to \Aut(\cs(G))$ denote the dual action of
$\hgab$ on $\cs(G)=\CC\rtimes_{\id}G$.
Then $(A\rtimes_{\alpha}G, \hgab, \widehat{\alpha})$
is isomorphic to the $C_0(X)$-system
$(A\otimes_X(Z\times_{\hgab}\cs(G)),\hgab,
\id\otimes_X\Ind\mu)$.
\end{thm}
\begin{proof} Choose a locally unitary action
$\delta:G\to\Aut(C_0(X,\K))$ such that $\zeta(\alpha)=\zeta(\delta)$.
Then it follows from Remarks
\ref{rem-stable}~and
\ref{rem-locuni} that $\alpha\otimes\id$ is exterior equivalent
to $\id\otimes_X\delta$ and that $\alpha$ is locally exterior equivalent
to $\id$. Thus, Theorem~\ref{thm-mostimportant} implies
that $(A\rtimes_{\alpha}G, \hgab, \widehat{\alpha})$ is
isomorphic to
$\(\fc Z*{(A\rtimes_{\id}G)},\hgab,\fc Z*{\widehat{\id}}\)$.
But
$(A\rtimes_{\id}G,\hgab, \widehat{\id})$ is $C_0(X)$-isomorphic to
$(A\otimes_{\max}\cs(G), \hgab, \id\otimes\mu)$, and the latter may
be written as
$(A\otimes_X(C_0(X,\cs(G)), \hgab, \id\otimes_X\mu)$.
But then we can apply Propositions \ref{prop-tensor}~and
\ref{proptriv} to see that
$(A\rtimes_{\alpha}G, \hgab, \widehat{\alpha})$ is
$C_0(X)$-isomorphic to
$(A\otimes_X(Z\times_{\hgab}\cs(G)),\hgab, \id\otimes_X\Ind\mu)$.
\end{proof}

We close this section with two corollaries about the dual
spaces of crossed products by locally unitary actions.
The first of these follows at once from \thmref{thm-loc-uni} together
with Remark~\ref{rem-final} and Proposition~\ref{prop-dualproduct}.

\begin{cor}\label{cor-locuniprim}
Let $A\in \CR(X)$ and let $\alpha:G\to \Aut(A)$ be locally unitary.
Let $q:Z\to X$ be a $\hgab$-principal bundle corresponding to
$\zeta(\alpha)$. If either $A$ or $\cs(G)$ is nuclear (in particular, if
$A$ is type I or
$G$ is amenable), then $\Prim(A\rtimes_{\alpha}G)$
is isomorphic to $\Prim(A)\times_X\(Z\times_{\hgab}\Prim(\cs(G))\)$
as a topological bundle over $X$ with group $\hgab$.
In particular, if $\Prim(A)$  is Hausdorff (so that $\Prim(A)=X$),
then $\Prim(A\rtimes_{\alpha}G)$ is isomorphic to
$Z\times_{\hgab}\Prim(\cs(G))$ as a topological bundle
over $X$ with group~$\hgab$.
\end{cor}

If $A$ is of type I, then
we get the following result.

\begin{cor}\label{cor-locunidual}
Let $(A,G,\alpha)$ and $q:Z\to X$ be as in Corollary~\ref{cor-locuniprim}.
If in addition, $A$ is of type I,
then $\specnp{(A\rtimes_{\alpha}G)}$
is isomorphic to $\hA\times_X(Z\times_{\hgab}\hat{G})$
as a topological bundle over $X$ with group $\hgab$.
In particular, if $\hA$ is Hausdorff (so that $\hA=X$),
then $\specnp{(A\rtimes_{\alpha}G)}$ is isomorphic to
$Z\times_{\hgab}\hat{G}$ as a topological bundle
over $X$ with group~$\hgab$.
\end{cor}

%% file: loc-inner.tex
%
%

\section{Crossed products by locally inner
actions of smooth groups}\label{loc-inner}

In addition to our separability proviso, we shall want some additional
assumptions and notations to be in effect throughout this section.

\begin{ga}
We assume that
$G$ is a smooth second countable locally
compact group such that $\gab$ is compactly generated.  We fix a
representation group
\[
e\arrow{e} C \arrow{e} H \arrow{e} G \arrow{e} e,
\]
and identify $C$ with $\specnp{H^2(G,\TT)}$.  Therefore, we may assume
that the transgression map
$\tg:\widehat{C}=H^2(G,\TT)\to H^2(G,\TT)$ is the identity map
\cite[Remark~4.4]{ew2}.
We let
$\sigma\in Z^2(G, \specnp{H^2(G,\TT)})$ be the corresponding cocycle
in Moore cohomology \eqref{eq-sigma1}.  When convenient, we will view
$\sigma$ as an element of $Z^2\(G,C(H^2(G,\T),\T)\)$.
Finally, $X$ will be a second countable locally compact space,
and $A\in \CR(X)$.
As in previous sections, we will identify $A\tensor\K$ and
$A\xtensor\coxk$.
\end{ga}

Under these assumptions we were able to extend some
results of Packer \cite{judymc} and give a
classification of the exterior equivalence classes $\LI_G(A)$ of locally
inner actions of
$G$ on $A$ in terms of the cohomology groups $H^1(X,\shgab)$
and
$H^2(G,\TT)$ \cite[Theorem~6.3]{ew2}.
As a special case, it is useful to note that it follows
from the argument in \cite[Corollary 2.2]{ros2} that if
 $A$ has continuous trace, then the locally inner actions
of $G$ on $A$ coincide with the $C_0(X)$-linear actions on $A$.
In particular, $\LI_G(C_0(X,\K))$
coincides with the abelian group $\E_G(X)$ of exterior equivalence classes
of $C_0(X)$-linear actions of $G$ on $C_0(X,\K)$.
The group operation in $\E_G(X)$ is given by
$[\gamma]\cdot[\delta]:=[\gamma \xtensor\delta]$
\cite[\S5]{ew2}.

In this section, we want
to use our classification of locally inner actions
to describe the $C_0(X)$-bundle
structures of the crossed products $A\rtimes_{\alpha}G$ in terms
of $\cs(H)$.
To do this we have to recall the basic ingredients of our
classification theory.

If $[\alpha]\in \LI_G(A)$, then $[\alpha]$
determines a  continuous map $\varphi_{\alpha}:X\to
H^2(G,\TT)$ such that
the action $\alpha^x:G\to \Aut(A_x)$
is implemented by an $\varphi_{\alpha}(x)$-homomorphism
$v^x:G\to \UM(A_x)$ (compare with Definition~\ref{def-u-homomorphism}).
Now let
$\sigma\in Z^2(G,C(H^2(G,\TT),\TT))$
be the cocycle as in our standing assumptions.
Then we can pull back
$\sigma$ to a cocycle $\varphi_{\alpha}^*(\sigma)\in Z^2(G,C(X,\TT))$
given by
\[
\varphi_{\alpha}^*(\sigma)(s,t)(x)=\sigma(s,t)(\varphi_{\alpha}(x)).
\]
By \cite[Proposition 3.1]{horr}, there exists an
inner action $\gamma:G\to \Aut(C_0(X,\K))$ which is
implemented by a $\varphi_{\alpha}^*(\sigma)$-homomorphism
$v:G\to \UM(C_0(X,\K))$. If $[\gamma^o]$ denotes the inverse
of $[\gamma]$ in $\E_G(X)$, then $\alpha\otimes_X\gamma^o$
is locally unitary and, therefore,
we obtain a class $\zeta_H(\alpha):=\zeta(\alpha\otimes_X\gamma^o)\in
H^1(X,\shgab)$ which determines the exterior equivalence class of
$\alpha\otimes_X\gamma^o$
(compare with the discussion in the previous section).
It follows from \cite[Theorem~6.3]{ew2} that $[\alpha]
\mapsto \zeta_H(\alpha) \oplus \varphi_{\alpha}$ defines an injection
\[
\Phi^H:\LI_G(A)\to H^1(X,\shgab)\oplus C(X,H^2(G,\TT)).
\]
Furthermore, $\Phi$ is a bijection whenever $A$ is stable.
(When $A=C_0(X,\K)$, the map $\Phi^H:\E_G(X)\to
H^1(X,\shgab)\oplus C(X,H^2(G,\TT))$ is an isomorphism of abelian
groups \cite[Theorem~5.4]{ew2}.)

\begin{definition}
[{\cite[Theorem~5.4 and Lemma~6.1]{ew2}}]
\label{def-zetaphi}
Let $\varphi_{\alpha}\in C(X, H^2(G,\TT))$  and $\zeta_H(\alpha)\in
H^1(X,\shgab)$ be
as above. Then  we say that $\varphi_{\alpha}:X\to H^2(G,\TT)$
is the {\em Mackey obstruction map\/} for $\alpha$ and we say
that $\zeta_H(\alpha)$ is the
{\em Phillips-Raeburn obstruction of $\alpha$
with respect to $H$}.
\end{definition}

\begin{remark}\label{rem-reprgroup}
(a) Notice that the class of the cocycle $\sigma$ in $H^2(G,
C(H^2(G,\TT),\TT))$ in our standing assumptions depends on
the choice of the representation group $H$.
This implies that in general, the action $\gamma$ constructed above, and
therefore the class
$\zeta_H(\alpha)\in H^1(X,\shgab)$, depends on the choice of $H$.
However, the class of $\sigma$ is uniquely
determined by the choice of $H$, which implies that the class of
$\varphi_{\alpha}^*(\tilde{\sigma})$ in $H^2(G, C(X,\TT))$ and hence
the exterior equivalence class of $\gamma$ is also uniquely
determined by $H$ (compare with Remark~\ref{rem-inner}).

(b) If $\delta:G\to \Aut(C_0(X,\K))$ is a locally unitary action
with $\zeta(\delta)=\zeta_H(\alpha):=\zeta(\alpha\otimes_X\gamma^o)$,
then it follows from the above constructions and the
classification of locally inner actions via $H^1(X,\shgab)$, that
$\alpha\otimes\gamma^o$ is exterior equivalent to $\id_A\otimes_X\delta$
as actions on $A\otimes \K$. Further, since $[\gamma^o]$ is the inverse
of $[\gamma]$ in $\E_G(X)$, it follows that $\gamma^o\otimes_X\gamma$
is exterior equivalent to the trivial action on $C_0(X,\K)$.
Hence we see that $\alpha\otimes\id_{\K}$ is exterior
equivalent to both
$\id_A\otimes_X\delta\otimes_X\gamma$ and $
\id_A\otimes_X\gamma\otimes_X\delta$.
\end{remark}

The main idea in our description of crossed products by locally inner
actions is to use the group $\cs$-algebra $\cs(H)$ as a universal bundle
over the locally compact space $H^2(G,\TT)$ with fibres isomorphic to the
twisted group algebras $\cs(G,\om)$, $[\om]\in H^2(G,\TT)$.
The fact that a representation group $H$ of $G$
does provide a bundle over
$H^2(G,\TT)$ with fibres $\cs(G,\om)$ was first observed by
Packer and Raeburn in \cite[\S1]{para3}.

More generally, let $e\to N\to L\to G\to e$
be any locally compact central extension of $G$ by an abelian group
$N$. Let $i_L:L\to \UM(\cs(L))$ denote the canonical map
and let $\phi:\cs(N)\to \M(\cs(L))$ denote the
integrated form of the restriction of $i_L$ to $N$.
Since, by assumption, $N$ is central in $L$ it follows that
$\phi$ takes image in $\ZM(\cs(L))$. Thus, if we
identify $\cs(N)$ with $C_0(\widehat{N})$ via the Gelfand transform,
we see that $\cs(L)$ has a canonical structure as a
$C_0(\widehat{N})$-algebra. In order to see that the fibres
are exactly what we want,
it is convenient to write $\cs(L)$ as a central twisted crossed product as
in the following lemma. Note that for any abelian locally compact group
$N$, we may view $N$ as a closed subgroup of
$C(\widehat{N},\TT)$ by identifying $N$ with the Pontryagin
dual of $\widehat N$.

\begin{lem}\label{lem-centralcross}
Let $e\to N\to L\to G\to e$ be a second countable central extension
of $G$ by the abelian group $N$, and let
$\sigma\in Z^2(G,N)\subseteq Z^2(G,C(\widehat{N},\TT))$ be given by
\[
\sigma(s,t)(\chi)=\chi(c(s)c(t)c(st)^{-1})
\]
for some Borel section $c:G\to L$ satisfying  $c(eN)=e$.
Let $\phi:C_0(\widehat{N})\to \ZM(\cs(L))$ denote the
canonical map described above and let $v:G\to \UM(\cs(L))$ be
given by $v(s)=i_L(c(s))$. Then the following assertions are true:
\begin{enumerate}
\item
$(\phi, v)$ is a covariant homomorphism of the twisted system
$(C_0(\widehat{N}), G,\id, \sigma)$ whose integrated
form $\phi\rtimes v$ is a $C_0(\widehat{N})$-isomorphism from
$\cs(G, \widehat{N}, \sigma)$ onto $\cs(L)$.
\item
For each $\chi\in \widehat{N}$, the fibre $\cs(L)_{\chi}$ is isomorphic
to the twisted group algebra $\cs(G, \tg(\chi))$.
\item
If  $\mu:\gab\to \Aut(\cs(L))$ is given via restriction of
the dual action of $\hlab$ to the closed subgroup
$\hgab$ of $\hlab$, then
$\phi\rtimes v$ intertwines $\widehat{\sigma}$ and $\mu$.
\item
If $G$ is amenable, then $\cs(L)$ is a continuous
$C_0(\widehat{N})$-bundle.
\end{enumerate}
\end{lem}

Part~(a) of
the above result is a very special case of \cite[Theorem 4.1]{para1},
and can also be deduced from \cite[Proposition 5.1]{para1}
using the decomposition of group algebras by Green's twisted
crossed products. Part~(b) is a special case of
\cite[Theorem~1.2]{para2}. Since this result and the special form of
the isomorphism will be crucial here, we give a direct proof for
convenience.

\begin{proof}[Proof of Lemma~\ref{lem-centralcross}]
To prove~(a), it will suffice to show that $(\cs(L), \phi, v)$ is a crossed
product for the twisted system $(C_0(\widehat{N}), G,\id, \sigma)$
(Definition~\ref{def-twistcross}).  Recall that we have identified
$C_0(\widehat{N})$ with
$\cs(N)$ via Gelfand  transform.  This transforms $\sigma: G\times
G\to N\subseteq C(\widehat{N},\TT)=\UM(C_0(\widehat{N}))$ to the cocycle
$\tilde{\sigma}:G\times G\to \UM(\cs(N))$ given by
$\tilde{\sigma}(s,t)=i_N(\sigma(s,t))$, where $i_N:N\to\UM(\cs(N))$
denotes the natural map.
Thus we have to check that $(\cs(L), \phi, v)$ is a crossed
product for $(\cs(N), G,\id,\tilde{\sigma})$, where
$\phi:\cs(N)\to \ZM(\cs(L))$ denotes the integrated form of $i_H\restr{N}$.
It follows from
$$v(s)v(t)=i_L(c(s)c(t))=i_L(c(s)c(t)c(st)^{-1})i_L(c(st))=
\phi(\tilde\sigma(s,t))v(st)$$
that
$(\phi,v)$ is a covariant homomorphism
of $(\cs(N), G, \id, \tilde{\sigma})$ into $\M(\cs(L))$;
this implies the first condition of Definition
\ref{def-twistcross}. In order to check the second, choose Haar
measures on $G$, $L$, and $N$ so that we may view $L^1(G\times N)$ as
a dense subset of $L^1(G,\cs(N))$.
Then, if $f\in L^1(G\times N)$,
$$
\phi\rtimes v(f)=\int_G \phi(f(s,\cdot))v(s)\,ds
=\int_G \int_N f(s,n) i_L(n)i_L(c(s))\,dn\,ds
=\int_L \tilde{f}(l) i_L(l)\,dl,
$$
where we define $\tilde{f}\in L^1(L)$ by $\tilde{f}(c(s)n)=f(s,n)$.
Since $(s,n)\mapsto c(s)n: G\times N\to L$ is a Borel isomorphism,
it follows that $\phi\rtimes v\big(L^1(G, \cs(N))\big)$
is a dense subspace of $\cs(L)$.
Finally, we have to check that any
covariant homomorphism $(\Psi, w)$ of $(\cs(N), G, \id, \tilde{\sigma})$
into $\M(D)$, for some separable $\cs$-algebra $D$,
has an integrated form $\Psi\rtimes w:\cs(L)\to M(D)$ satisfying
$(\Psi\rtimes w)\circ \phi=\Psi$ and $(\Psi\rtimes w)\circ v=w$.
For this let $\Psi(n):=\Psi(i_N(n))$ for $n\in N$.
Then define $u:L\to\UM(D)$ by
\[
u(c(s)n)=w(s)\Psi(n)
\quad
\text{for all $s\in G$, $n\in N$.}
\]
Using the covariance condition
for $(\Psi, w)$, a computation shows that $u$ is a strictly Borel
homomorphism,
and that the integrated form of $u$ has the properties
required of $\Psi\rtimes w$.

By the uniqueness of the crossed product it follows that
$\phi\rtimes v:\cs(G, \widehat{N},{\sigma})\to \cs(L)$
is an isomorphism. Since
$\phi:C_0(\widehat{N})\to \M(\cs(L))$ is clearly
$C_0(\widehat{N})$-linear,
it follows from Lemma~\ref{lem-CXcov} that $\phi\rtimes v$ is
$C_0(\widehat{N})$-linear. This proves~(a).

Since evaluation of $\tg(\chi)$ is defined via evaluation of $\sigma$
at $\chi$, part~(b) of the Lemma follows from \cite[Theorem 5.1]{may1},
and~(c) follows from
the definitions of $v$ and the dual actions.
The final assertion
follows from \cite[Theorem 1.2]{para2}.
\end{proof}

Since the representation group is a central extension of $G$ by $C$
and since $\widehat C$ has been identified with $H^2(G,\T)$, the
following is an immediate corollary.

\begin{cor}[{cf., \cite[Section 1]{para2}}]\label{cor-rep}
Let $G$, $H$ and $\sigma$ be as in our standing assumptions.
Then $\cs(H)$ is a $C_0\(H^2(G,\TT)\)$-algebra
which is $C_0\(H^2(G,\TT)\)$-linearly and $\hgab$-equivariantly
isomorphic to $\cs(G, H^2(G,\TT),\sigma)$. The fibres
$\cs(H)_{[\om]}$ are isomorphic to $\cs(G,\om)$ for each
$[\om]\in H^2(G,\TT)$, and, if $G$ is amenable, then $\cs(H)$ is a
continuous $C_0\(H^2(G,\TT)\)$-bundle.
\end{cor}

Recall that if $A$ is a $C_0(X)$-algebra, $B$ is a $C_0(Y)$-algebra
and $f:X\to Y$  is a continuous map, then $A$ becomes a $C_0(Y)$-algebra
via composition with $f$.
Thus we can form the balanced tensor product
$A\otimes_fB:=A\otimes_YB$
which becomes a $C_0(X)$-algebra via composition with
the natural map $i_A:A\to \M(A\otimes_fB)$.
Therefore, for any $C_0(X)$-algebra $A$ we have
\begin{equation}\label{eq-pull-back}
A\otimes_Xf^*B=A\otimes_X(C_0(X)\otimes_fB)\cong(A\otimes_XC_0(X))\otimes_fB
\cong A\otimes_fB,
\end{equation}
where $f^*B$ is the usual pull-back $C_0(X)\tensor_f B$
(Remark~\ref{rem-final}).
Moreover, if $\beta:G\to \Aut(B)$ is a $C_0(Y)$-linear action
of a group $G$ on $B$, then we obtain a $C_0(X)$-linear action
$\id\otimes_f\beta$ (resp., $f^*\beta$) of $G$ on $A\otimes_fB$
(resp., $f^*B$), and the isomorphisms in equation (\ref{eq-pull-back})
are $C_0(X)$-linear and $G$-equivariant.

The following lemma shows that the pull-back of a twisted transformation
group \cs-algebra, is itself a twisted transformation group algebra twisted
by the pull-back of the original
cocycle:
if $u\in Z^2(G,C(Y,\T))$ and $f:X\to Y$, then $f^*(u)\in
Z^2(G,C(X,\T))$ is defined by $f^*(u)(s,t)(x)=u(s,t)\(f(x)\)$.

\begin{lem}\label{lem-pull-back}
Let $f:X\to Y$ be a continuous map between the second countable
locally compact spaces $X$ and $Y$, and let $u\in Z^2(G, C(Y,\TT))$.
Then there exists a $C_0(X)$-linear and $\hgab$-equivariant isomorphism
between
$\cs(G,X,f^*(u))$ and $f^*(\cs(G,Y,u))$.
\end{lem}
\begin{proof}
Let $\Phi:C_0(X)\otimes_YC_0(Y)\to C_0(X)$ denote the isomorphism given
on elementary tensors by $\Phi(h\otimes_Yg)(x)=h(x)g(f(x))$.
Then
$$\Phi(1\otimes_Yu(s,t))(x)=u(s,t)(f(x))=f^*(u)(s,t),$$
which implies that $\Phi$ transforms the twisted action
$(\id\otimes_Y\id, 1\otimes u)$ to $(\id, f^*(u))$.
Thus Proposition~\ref{prop-twisttensor} implies that
\begin{align*}
\mathcs(G,X,f^*(u))&=C_0(X)\rtimes_{\id,f^*(u)}G\cong
(C_0(X)\otimes_YC_0(Y))\rtimes_{\id\otimes_Y\id, 1\otimes_Yu}G\\
&\cong C_0(X)\otimes_Y(C_0(Y)\rtimes_{\id,u}G)
= f^*(\mathcs(G,Y,u)).
\end{align*}
Each of the above isomorphisms is $C_0(X)$-linear and $\hgab$-equivariant;
the first because $\Phi$ is clearly $C_0(X)$-linear and the second
due to Proposition~\ref{prop-twisttensor}.
\end{proof}

We are now prepared to state our main result.

\begin{thm}\label{thm-main}
Let $A$, $G$ and $H$ be as in our standing assumptions and let
$\alpha:G\to \Aut(A)$ be a locally inner action. Let
$\varphi_{\alpha}\in C(X,H^2(G,\TT))$ be the Mackey obstruction map for
$\alpha$ and let
$\zeta_H(\alpha)\in H^1(X,\shgab)$ be the cohomology class corresponding to
$\alpha$ with respect to $H$.
Further, let $q:Z\to X$ be a principal $\hgab$-bundle corresponding to
$\zeta(\alpha)$. Then
there exists a $C_0(X)$-linear and
$\hgab$-equivariant isomorphism between
$A\rtimes_{\alpha}G$ and $\fc Z*{(A\otimes_f\cs(H))}$, where
where $f:X\to H^2(G,\TT)$ is
defined by $f(x)=\varphi_{\alpha}(x)^{-1}$ for all $x\in X$.
\end{thm}

\begin{proof}
Let $\sigma\in Z^2\(G,C(H^2(G,\TT),\T)\)$ be as in our standing
assumptions,
and let $\varphi_{\alpha}^*(\sigma)\in Z^2(G, C(X,\TT))$ denote the
pull-back of $\sigma$ via $\varphi_{\alpha}$. Choose
a $\varphi_{\alpha}^*(\sigma)$-homomorphism $v:G\to \UM(C_0(X,\K))$
and let $\gamma:G\to \Aut(C_0(X,\K))$ denote the inner action
implented by $v$.  Further let
 $\delta:G\to \Aut(C_0(X,\K))$ be a locally unitary
action such that $\zeta(\delta)=\zeta(\alpha)$.
Then it follows from Remark~\ref{rem-reprgroup} that
$\alpha\otimes\id_{\K}$ and
$\id_A\otimes_X\gamma\otimes_X\delta$ are exterior equivalent as
actions on
$A\otimes\K$. The discussion following
Definition~\ref{def-u-homomorphism}  implies that $v$
implements an exterior equivalence between $\gamma=(\gamma,1)$ and the
twisted action
$\big(\id_{\coxk},\overline{\varphi_\alpha^*(\sigma})\big)$.   If we
identify $\coxk$ and $\cox\xtensor\coxk$ and view $f^*(\sigma)$ as taking
values in $\UM\(\cox\)$, then $\gamma$ is exterior equivalent to
$(\id_{\cox}\xtensor\id_{\coxk},f^*(\sigma)\xtensor 1)$. If we let $\tilde
f^*(\sigma)$ be the cocyle taking values in
$\UM(A)$ corresponding to $1\xtensor f^*(\sigma)$ via the
identification of $A$ with $A\xtensor\cox$, then we conclude that the
action $\alpha\tensor\id_\K$ is exterior equivalent to the twisted
action $\(\id_A\xtensor \delta,\tilde f^*(\sigma)\)$

We claim that $\alpha$ is locally exterior equivalent to
the twisted action $(\id_A, f^*(\sigma))$. To see this let
$x\in X$ and choose an open neighborhood $W$ of $x$ such
that $\alpha^W:G\to\Aut(A_W)$ is inner.
Then there exists a cocycle $u\in Z^2(G,C(W,\TT))$ and
a $u$-homomorphism $w:G\to\UM(A_W)$ such that $\alpha^W=\Ad w$,
and $w$ implements an exterior equivalence between
$\alpha^W$ and $(\id_A,\bar{u})$. We want to show that
there exists a possibly smaller neighborhood $W_1$ of
$x$ such that the restrictions of $\bar{u}$ and $f^*(\sigma)$ to
$W_1$ are cohomologous, or, equivalently, such that
the product $u\cdot f^*(\sigma)$ restricted to $W_1$
is cohomologous to the trivial cocycle. This would imply the claim
since the twisted actions $(\id_A,\bar{u})$ and  $(\id_A,
f^*(\sigma))$ would then be exterior equivalent when restricted to
$A_{W_1}$.
Since evaluation of $[u(y)]\in H^2(G,\TT)$ of $u$
at a given point $y\in W$ must coincide with $\varphi_{\alpha}(y)$
(since both correspond to the same inner action on the fibre $A_y$),
and since $[f^*(\sigma)(y)]=\varphi_{\alpha}(y)^{-1}$, it follows
that $(u\cdot f^*(\sigma))(y)$ is cohomologous to the trivial cocycle
for every $y\in W$, hence $u\cdot f^*(\sigma)$ is pointwise trivial.
Thus, by Rosenberg's theorem \cite[Theorem 2.1]{ros2}, it follows
that $u\cdot f^*(\sigma)$ is locally trivial, which is precisely
what we want.

Since $\alpha\otimes\id$ is exterior equivalent to
$(\id_A\otimes_X\delta,\tilde f^*(\sigma)\otimes 1)$ we can apply
Theorem~\ref{thm-mostimportant} to obtain a
$C_0(X)$-linear and $\hgab$-equivariant isomorphism between
$A\rtimes_{\alpha}G$ and $\fc Z*{(A\rtimes_{\id,\tilde f^*(\sigma)}G)}$.
But $A\rtimes_{\id,\tilde f^*(\sigma)}G$ is $C_0(X)$-linearly and
$\hgab$-equivariantly isomorphic to
$A\otimes_X\cs(G,X,f^*(\sigma))$ by Proposition~\ref{prop-twisttensor}
(since $A\rtimes_{\id,\tilde f^*(\sigma)}G \cong \big(\cox\xtensor A\big)
\rtimes_{\id\xtensor\id,f^*(\sigma)\xtensor 1}G$).
Finally, it follows from Corollary~\ref{cor-rep} and
Lemma~\ref{lem-pull-back} that
$\cs(G,X,f^*(\sigma))$ is $C_0(X)$-linearly and
$\hgab$-equivariantly isomorphic to $f^*(\cs(H))$.
Since
$A\otimes_f\cs(H)=A\otimes_Xf^*(\cs(H))$, this completes the proof.
\end{proof}

\begin{remark}\label{rem-main}
(a) Note that our theorem applies to all $C_0(X)$-linear
actions of a smooth group $G$ (with $\gab$ compactly generated)
on a separable continuous trace algebra
$A$ with spectrum~$X$.

(b) Since $A\otimes_f\cs(H)$ is isomorphic to
$A\otimes_Xf^*(\cs(H))$, since
$\fc Z*{(A\otimes_Xf^*(\cs(H))}$ is isomorphic to $A\otimes_X\big(\fc
Z*{f^*(\cs(H))}\big)$ by Proposition~\ref{prop-tensor}, and since
all the isomorphisms are $C_0(X)$-linear and $\hgab$-equivariant, we see
that the crossed product $A\rtimes_{\alpha}G$ is obtained  via the
iteration of the following basic bundle operations:
First take the pull-back $f^*(\cs(H))$ of the universal bundle $\cs(H)$
via the continuous map $f=\overline{\varphi_{\alpha}}:G\to H^2(G,\TT)$.
Then construct the $G$-fibre product $\fc Z*{f^*(\cs(H)})$.
Finally take the
fibre product (i.e., the balanced tensor product) of
$\fc Z*{f^*(\cs(H)})$ with $A$.

(c) If $A$ is type I, then this bundle-theoretic description of
$A\rtimes_{\alpha}G$ gives a bundle theoretic
description of the spectrum $\specnp{(A\rtimes_{\alpha}G)}$:
it is isomorphic (as a topological bundle with group $\hgab$)
to $\widehat{A}\times_X(Z*(f^* \widehat{H}))$, and if
$\widehat{A}=X$ we get
$\specnp{(A\rtimes_{\alpha}G)}=Z*(f^* \widehat{H})$.
If $A$ is nuclear, we obtain a similar description of
$\Prim(A\rtimes_{\alpha}G)$.
\end{remark}

We close this section with some corollaries of
Theorem~\ref{thm-main} for certain special cases.
For example, combining \thmref{thm-main} with \lemref{lem-trivial}
immediately yields the following corollary.  Note that $H^1(X,\shgab)$
is always trivial if $\gab$ is a vector
group, or if $X$ is contractible (cf., e.g.,
\cite[Corollary~4.10.3]{husemoller}).

\begin{cor}\label{cor-trivialh1}
Suppose that $H^1(X,\shgab)$ is trivial, and let $A$,
$G$ and
$H$ be as in our general assumptions. Then
$A\rtimes_{\alpha}G$ is isomorphic to $A\otimes_X f^*(\cs(H))$
for any locally inner action $\alpha:G\to \Aut(A)$,
where $f(x)=\varphi_{\alpha}(x)^{-1}$ for $x\in X$.
\end{cor}

So in the case where $H^1(X,\shgab)$ is trivial we get a description
of $A\rtimes_{\alpha}G$ by pulling back the universal bundle
$\cs(H)$ via $f$ and then taking the fibre product with $A$.

Another interesting special situation occurs when
$\varphi_{\alpha}$ is constant. If $A$ has continuous trace, such
systems were called  {\em pointwise projective unitary\/} in
\cite{echros}.
In this situation the class $\zeta(\alpha)$ does
not depend on  the choice of the representation group $H$; if
$[\om]\in H^2(G,\TT)$ is the constant value of $\varphi_{\alpha}$,
then the
pull-back of $\sigma$ via $\varphi_{\alpha}$ would always give a cocycle
cohomologous to the cocycle $\varphi_{\alpha}^*(\om)$ defined by
$\varphi_{\alpha}^*(\om)(s,t)(x)=\om(s,t)$.

\begin{cor}\label{cor-proloc}
Let $G$ be smooth, $A\in \CR(X)$, and let $\alpha:G\to \Aut(A)$ be a
locally inner action of $G$ on $A$ so that there is a class $[\om]\in
H^2(G,\TT)$ with
$\varphi_{\alpha}(x)=[\om]$ for all $x\in X$. Let
$q:Z\to X$ be a principal bundle corresponding to $\zeta(\alpha)$.
Then $A\rtimes_{\alpha}G$ is $C_0(X)$-linearly and
$\hgab$-equivariantly isomorphic to
$A\otimes_X\(Z\times_{\hgab}\cs(G,\bar{\om})\)$.
Moreover, if $A$ is type I, then $\specnp{(A\rtimes_{\alpha}G)}$
is isomorphic to $\widehat{A}\times_X
(Z\times_{\hgab}\widehat{G}_{\bar{\om}})$,
where $\widehat{G}_{\bar{\om}}$ denotes the set of equivalence
classes of irreducible $\bar{\om}$-representations of $G$.
\end{cor}

\begin{proof} If $H$ is any representation group of $G$
as in our standing assumptions, then since
$f(x)=\varphi_{\alpha}(x)^{-1}$, it is easily seen that
$f^*(\cs(H))=C_0(X, \cs(G,\bar{\om}))$.
Thus the result follows from
Theorem~\ref{thm-main} and Proposition~\ref{prop-tensor}.
\end{proof}

Note that the above result also holds true under the weaker
assumption that $\alpha$ is {\em locally projective unitary}
in the sense that there exists an action $\beta:G\to\Aut(\K)$
with Mackey obstruction $[\bar{\om}]\in H^2(G,\TT)$ such that
$\alpha\otimes\beta:G\to \Aut(A\otimes\K)$ is locally unitary
(compare with \cite[p. 38]{echros}). The proof is similar to the proof
of Theorem~\ref{thm-loc-uni} and we omit the details.

%% file: action.tex
%
%
\section{Some applications to actions of $\RR^n$}\label{action}

As an application of our results, we want to consider locally inner
actions of $\R^n$ on \cs-algebras $A$ with Hausdorff
spectrum~$X$.  (Since $A$ is assumed to be separable, it is
necessarily type~I if its spectrum is Hausdorff.)
When $A$ has continuous-trace, this is equivalent to
requiring that $(A,\R^n,\alpha)$ is a \cox-system.
A representation group $H$ for $\RR^n$ was explicitly
constructed in \cite[Example 4.7]{ew2}. $H$ is a simply connected and
connected two-step nilpotent Lie group (in fact it is the universal
two-step nilpotent group with $n$ generators), and by \cite[Proposition
4.8]{ew2}
it is unique up to isomorphism.
Notice that in case $n=2$, $H$ is just the real Heisenberg group of dimension
three. Since $H^1(X,{\widehat{\R}}^n)$ is trivial for all second countable
locally compact spaces $X$ the following is a consequence of
\corref{cor-trivialh1}.

\begin{thm}\label{thm-rn}
Let $A$ be a separable type I $\cs$-algebra with Hausdorff spectrum $X$
and let $\alpha:\RR^n\to \Aut(A)$ be a locally inner action of
$\RR^n$ on $A$. Let $\varphi_{\alpha}:X\to H^2(\RR^n,\TT)$ denote
the Mackey obstruction map and let $f(x)=\varphi_{\alpha}(x)^{-1}$ for
$x\in X$.
Then $A\rtimes_{\alpha}\RR^n$ is $C_0(X)$-linearly and
${\widehat{\R}}^n$-equivariantly isomorphic to $A\otimes_Xf^*(\cs(H))$
and $\specnp{(A\rtimes_{\alpha}\RR^n)}$ is isomorphic to
$f^*(\hat{H})$ as a topological bundle over $X$ with group
${\widehat{\R}}^n$.
\end{thm}

Notice that it is well known that $A\rtimes_{\alpha}\RR^n$ is type I
in the situation above. This can also be deduced from the above result,
since any connected nilpotent Lie group $N$ has type I group $\cs$-algebra
$\cs(N)$.
We want to use our result to obtain a more detailed description
of $A\rtimes_{\alpha}\RR^n$ and its spectrum.
Recall that if $G$ is an abelian group and $\om\in Z^2(G,\TT)$,
then the {\em symmetry group} $\Sigma_{\om}$ of $\om$ is defined
as  $\Sigma_{\om}:=\{s\in G: \om(s,t)=\om(t,s)\;\text{for all}\; t\in G\}$.
Note that the symmetry groups only depend on the cohomology classes
$[\om]\in H^2(G,\TT)$. If $G=\RR^n$ then any cocycle of
$G$ is cohomologous to a cocycle of the form $\om(s,t)=e^{iJ(s,t)}$,
where $J$ is a skew symmetric form on $\RR^n$
\cite[Theorem~10.38]{var}.
It follows that the symmetry group equals the radical of $J$, so that
$\Sigma_{\om}$ is actually a vector subgroup of $\RR^n$.  The symmetry
groups play an important r\^ole in the representation theory of
two-step nilpotent groups and crossed products by abelian groups (see
for instance \cite{bagpac, echros, ech4}). For example, if
$A$ is a separable (type I) $\cs$-algebra with Hausdorff spectrum
$\hA=X$ and if $\alpha:G\to \Aut(A)$ is a $C_0(X)$-linear action of
the second countable abelian group $G$ on $A$, then each fibre
$\Prim(A_x\rtimes_{\alpha^x}G)$ of
$\Prim(A\rtimes_{\alpha}G)$ over $x\in X$ is
$\widehat{G}$-equivariantly homeomorphic to $\widehat{\Sigma}_x$,
where $\Sigma_x$ denotes the symmetry group of the
Mackey-obstruction $\varphi_{\alpha}(x)$
(see \cite[Theorem 1.1]{horr}).

Thus, if we view $\Prim(A\rtimes_{\alpha}G)$
as a topological bundle over $X$ with group $\hat{G}$,
then there is a nice description of the fibres.
The problem is to get the global picture of the bundle.
When $G=\RR^n$, we will deduce a description
from the following fact, which was observed by Baggett and Packer
in \cite[Remark 2.5]{bagpac}. We include a relatively short
proof using Kirillov
theory (see \cite{corgr}).

\begin{prop}[{cf, \cite{bagpac}}]\label{two-step-dual}
Let $H$ be a connected and simply connected two-step nilpotent Lie group
with center $Z$ and quotient $G=H/Z$. Let $\tg:\widehat{Z}\to H^2(G,\TT)$
denote the transgression map and let $\Sigma_{\chi}$ denote the symmetry
group of
$\tg(\chi)$ for all $\chi\in \widehat{Z}$.
Then the topological bundle  $\widehat{H}$ over $\widehat{Z}$ with group
$\widehat{G}$  is isomorphic to the quotient space
$(\widehat{G}\times\widehat{Z})/\!\!\sim$, where
$\sim$ denotes the equivalence relation
$$(\mu,\chi)\sim (\mu',\chi')\Leftrightarrow \chi=\chi'\;\text{and}\;\;
\bar{\mu}\mu'\in \Sigma_{\chi}^{\perp},$$
and where $(\widehat{G}\times\widehat{Z})/\!\!\sim$ is equipped with the
canonical
structure as a topological bundle over $X$ with group $\widehat{G}$.
\end{prop}
\begin{proof}
Let $\mathfrak h$ denote the Lie algebra of $H$. Since $H$ is two step
nilpotent, we may identify $H$ with $\mathfrak h$, where multiplication
on $\mathfrak h$ is given by the Campbell-Baker-Hausdorff formula
$X\cdot Y=X+Y+\frac{1}{2}[X,Y]$.  We also identify $Z$ with the center
$\mathfrak z$ of $\mathfrak h$ and $G$ with $\mathfrak h/\mathfrak z$,
with addition as group operation.  We can write the vector space dual
of $\mathfrak h$ as $\mathfrak h^*=\mathfrak z^*\oplus \mathfrak
z^{\perp}=\mathfrak z^*\oplus\mathfrak g^*$, and we may identify
$\widehat{Z}$ and $\widehat{G}$ with $\mathfrak z^*$ and $\mathfrak
g^*$ via the exponential map.  By Kirillov's theory we have
$\widehat{H}\cong \mathfrak h^*/\Ad^*(H)$, where the coadjoint action
$\Ad^*$ is
given in the two-step case by $\Ad^*(\exp(X))(f)=f- f\circ \ad(X)$ (recall that
$\ad(X)(Y)=[X,Y]$).
Hence the
$\Ad^*(H)$ orbit of
$f$ coincides with
$f+\mathfrak r_f^{\perp}$, where $\mathfrak r_f$ denotes the radical
of the skew symmetric form $B_f(X,Y)=f([X,Y])$ on $\mathfrak h$ (see
for instance \cite[Proposition 1.6]{lipros}).  Notice that $\mathfrak
r_f$ always contains $\mathfrak z$ and only depends on
$f\restr{\mathfrak z}$. So if we put $\mathfrak g_f=\mathfrak
r_f/\mathfrak z$ it follows that $\mathfrak h^*/\Ad^*(H)$ coincides
with the space $(\mathfrak g^*\times \mathfrak z^*)/\!\!\sim$, where
$$(g,f)\sim (g',f')\Longleftrightarrow f=f'\quad\text{and}\quad  g-g'\in
\mathfrak g_f^{\perp}.$$
Using the identifications  $\mathfrak g^*\cong\widehat{G}$
and $\mathfrak z^*\cong\widehat{Z}$, the result will follow
if we can show that the algebras $\mathfrak g_f$ (viewed as subgroups of $G$)
coincide with the symmetry groups.

For this let $X_1,\ldots, X_m, Y_1,\ldots , Y_n$ be a Jordan-H\"older
basis
for the Lie algebra $\mathfrak h$ such that
$\mathfrak g=\langle Y_1+\mathfrak z,\ldots, Y_n+\mathfrak z\rangle$.
Define a continuous cross
section
$c:\mathfrak g\to \mathfrak h$ by $c\big(\sum_{i=1}^na_i(Y_i+\mathfrak
z)\big) =
 \sum_{i=1}^na_iY_i$.
Then, for $Y=\sum_{i=1}^na_iY_i$ and  $Y'=\sum_{i=1}^nb_iY_i$,
 it follows from the Campbell-Baker-Hausdorff formula that
\[
c(Y+\mathfrak z)\cdot c(Y'+\mathfrak z)\cdot c(Y+Y'+\mathfrak z)^{-1}=
Y\cdot Y'\cdot(-Y-Y')=\frac{1}{2}[Y,Y'].
\]
Hence, if $f\in \mathfrak z^*$ is viewed as a character of $Z$, the cocycle
$\om_f\in Z^2(G,\TT)$ is given by $e^{if(\frac{1}{2}[Y,Y'])}$. Thus
the symmetry group coincides with the radical of the
skew symmetric form given by $f$ on $\mathfrak g$, which is precisely
$\mathfrak  g_f$.
\end{proof}

Combining this and Theorem~\ref{thm-rn} we get:

\begin{cor}\label{cor-rn}
Let $\alpha:\RR^n\to\Aut(A)$ be a locally inner action of $\RR^n$ on
a separable (type I) $\cs$-algebra with Hausdorff spectrum $\hA=X$.
For each $x\in X$ let $\Sigma_x$ denote the symmetry group of
$\varphi_{\alpha}(x)\in H^2(G,\TT)$. Then $\specnp{(A\rtimes_{\alpha}\RR^n)}$
is isomorphic to $({\widehat{\R}}^n\times X)/\!\!\sim$ as a topological
bundle over $X$ with group ${\widehat{\R}}^n$,
where $\sim$ is the equivalence relation
$$(\mu,x)\sim (\mu',x')\Longleftrightarrow x=x'\quad
\text{and}\quad \bar{\mu}\mu'\in
\Sigma_x^{\perp}.$$
\end{cor}
\begin{proof}
Let $H$ be the representation group of $\RR^n$. Since $H$ is a connected
and simply connected two-step nilpotent Lie group with center
$\specnp{H^2(\RR^n,\TT)}$, it follows from
Proposition~\ref{two-step-dual} that
$\widehat{H}$ is isomorphic (as a bundle) to $({\widehat{\R}}^n\times
H^2(\RR^n,\TT))/\!\!\sim$.
Moreover, by Theorem~\ref{thm-rn} we know that
$\specnp{(A\rtimes_{\alpha}\RR^n)}$ is isomorphic to
$f^*\widehat{H}\cong f^*(({\widehat{\R}}^n\times
H^2(\RR^n,\TT))/\!\!\sim\nolinebreak)$.
Since $f(x)=\varphi_{\alpha}(x)^{-1}$ and the symmetry groups of
$[\om]$ and $[\om]^{-1}$ coincide for all $[\om]\in H^2(G,\TT)$,
it is follows that $[\mu,x]\mapsto (x,[\mu,f(x)])$
is an isomorphism between $({\widehat{\R}}^n\times X)/\!\!\sim$ and
$f^*(({\widehat{\R}}^n\times H^2(\RR^n,\TT))/\!\!\sim)$.
\end{proof}

The previous result can fail for an arbitrary second countable
compactly generated
abelian group $G$; we know from the work of
Phillips-Raeburn and Rosenberg that if
$\alpha:G\to\Aut(A)$ is any $C_0(X)$-linear action of
$G$ on a continuous trace algebra $A$ with spectrum $X$
such that the Mackey obstruction map vanishes (i.e., $\alpha$
is pointwise unitary) that $\specnp{(A\rtimes_{\alpha}G)}$
can be any principal $\widehat{G}$-bundle, while
$(\widehat{G}\times X)/\!\!\sim$ is just the trivial bundle
$\widehat{G}\times X$ in this case.
On the other hand, it would be interesting to see whether the result
remains to be true if we replace $\widehat{G}\times X$ with an
appropriate principal $\widehat{G}$-bundle $q:Z\to X$.
That is, it would be interesting to know under what circumstances the
following question has a positive answer (see also
\cite[Remark~2.5]{bagpac}).

\begin{problem}
Let $\alpha:G\to\Aut(A)$ be a $C_0(X)$-linear action
of the second countable compactly generated abelian group $G$ on the
continuous trace algebra $A$ with spectrum $X$. Does there always
exist a principal $\widehat{G}$-bundle $q:Z\to X$ such that
$\Prim(A\rtimes_{\alpha}G)$ is isomorphic to $Z/\!\!\sim$ as a bundle
over $X$?
Here $\sim$ denotes the equivalence relation
$$z\sim z'\Longleftrightarrow q(z)=q(z')\quad\text{and} \quad
\bar{z}z'\in
\Sigma_x^{\perp},
$$
where $\bar{z}z'$ denotes the unique element $\chi$ of $\widehat{G}$ which
satisfies $\chi\cdot z= z'$.
\end{problem}

It is straightforward to check that
$Z/\!\!\sim$ is just the twisted bundle $Z*(\widehat{G}\times
X)/\!\!\sim$.
Clearly, the above problem is
strongly related to the problem of
describing the primitive ideal space of the
group
$\cs$-algebra of a two-step nilpotent group $H$ with center $Z$ and quotient
$H/Z=G$ as a quotient space of a principal $\widehat{G}$-bundle over
$\widehat{Z}$, as considered by Baggett and Packer in
\cite{bagpac}. On the one
hand, the problem for two-step nilpotent groups is a special case of
the above,
since by the Packer-Raeburn stabilization trick
\cite[Corollary~3.7]{para1}, we can
write $\cs(H)\otimes\K$ as a crossed product
$C_0(\widehat{Z},\K)\rtimes_{\beta}G$, for some $C_0(\widehat{Z})$-linear
action $\beta$. On the other hand, if the result is true for
the representation group $H$ of $G$, which is always two-step
nilpotent, then as in the proof of
Corollary~\ref{cor-rn}, one could
get the same result for all $C_0(X)$-linear actions of $G$
on separable continuous trace algebras with spectrum $X$
(or, more generally,
for locally inner actions on type I algebras with Hausdorff spectrum $X$).

We want to illustrate this for the special case $G=\ZZ^2$,
where the (unique)
representation group is the discrete Heisenberg group of rank
three.

\begin{thm}\label{thm-z2}
Let $\alpha:\ZZ^2\to \Aut(A)$ be a $C_0(X)$-linear action of $\ZZ^2$
on the separable continuous trace algebra $A$ with spectrum $X$.
Let $\zeta({\alpha})$ be the Phillips-Raeburn obstruction of $\alpha$
as defined in Definition~\ref{def-zetaphi}, and let
$q:Z\to X$ denote the corresponding principal $\TT^2=\widehat{\ZZ}^2$-bundle.
Then $\Prim(A\rtimes_{\alpha}\ZZ^2)$ is isomorphic to
$Z/\!\!\sim$ as a topological bundle over $X$ with group $\TT^2$.
\end{thm}
\begin{proof}
Recall that the discrete Heisenberg group $H$ is the set $\ZZ^3$
with multiplication given by $(n_1,m_1,l_1)(n_2,m_2,l_2)=
(n_1+n_2,m_1+m_2, l_1+l_2+n_1m_2)$.
The center $C$ of $H$ is given by $\{(0,0,l):l\in \ZZ\}$.
For each $t\in [0,1)$ let $\chi_t(l)=e^{i2\pi tl}$ denote the character
corresponding to $t$ under the identification of $\widehat{C}$ with
$\TT=\RR/\ZZ$.
Using the section $c:\ZZ^2\to H: c(n,m)=(n,m,0)$,
we easily compute that the cocycle $\om_t\in Z^2(\ZZ^2,\TT)$ corresponding
to $\chi_t\in \widehat{C}$ is given by
$$\om_t((n_1,m_1),(n_2,m_2))=e^{i2\pi tn_1m_2}.$$
If $t$ is irrational, then the symmetry group $\Sigma_t=\Sigma_{\om_t}$
is trivial, and if $t=\frac{p}{q}$, where $p$ and $q$ have no common
factors, then it is not hard to show that $\Sigma_t=q\ZZ\times q\ZZ\subseteq
\ZZ^2$. Thus it follows that $\om_t$ is identically $1$ when restricted
to the symmetry groups, and we can use \cite[Theorem 2.3]{bagpac}
to deduce that $\Prim(\cs(H))$ is isomorphic to
$(\widehat{\ZZ}^2\times\widehat{C})/\!\!\sim$ as a bundle over
$\widehat{C}$, where $\sim$ is the
usual equivalence relation.

By Theorem~\ref{thm-main} we have $A\rtimes_{\alpha}\ZZ^2=
Z*(f^*\cs(H))$, where $f(x)$ is the inverse of the Mackey obstruction
$[\om_x]$ for all $x\in X$. Hence
\begin{align*}
\Prim(A\rtimes_{\alpha}\ZZ^2)&=Z*\big(f^*\Prim(\cs(H))\big)
=Z*\big(f^*((\widehat{\Z}^2\times\widehat{C})/\!\!\sim))\big)\\
&=
Z*\big((\widehat{\Z}^2\times X)/\!\!\sim\big)= Z/\!\!\sim.
\qed
\end{align*}
\renewcommand{\qed}{}
\end{proof}

Finally, we point out that our results are
also helpful to the investigation of the structure
of continuous trace subquotients of the crossed
products $A\rtimes_{\alpha}\RR^n$, where $\alpha$ is a
$C_0(X)$-linear action on the continuous trace algebra $A$ with
spectrum $X$.
For this we first recall the following result due to the first author:

\begin{prop}[{\cite[Theorem 6.3.3]{ech4}}]
\label{cont-trace}
Let $A$ be a separable continuous-trace algebra with spectrum
$\hA=X$ and let $\alpha:\RR^n\to \Aut(A)$ be a $C_0(X)$-linear
action. Further, let $\dim:X\to \ZZ^+$ be defined by letting
$\dim(x)$ be the vector space dimension of $\Sigma_{x}$.
Then
$A\rtimes_{\alpha}\RR^n$ has continuous trace if and only if
$\dim:X\to \ZZ^+$ is continuous.
\end{prop}

More generally, if $\alpha:\RR^n\to \Aut(A)$ is any $C_0(X)$-linear
action of $\RR^n$ on a continuous trace algebra $A$ with spectrum $X$,
then there exists a finite decomposition series of ideals
$$\{0\}=I_0\subseteq I_i\subseteq\cdots \subseteq I_l=A
\rtimes_{\alpha}\RR^n,$$
with $l\leq \frac{n}{2}+1$ and all subquotients $I_k/I_{k-1}$, $1\leq
k\leq l$,
given by crossed products with continuous trace as in the
proposition above
(\cite[Theorem 6.3.3]{ech4}). We are now going to use our results
together with a recent result of Lipsman and Rosenberg
to compute explicitly the Dixmier-Douady invariant of these subquotients.

\begin{thm}\label{Dix-Douady}
Let $\alpha:\RR^n\to \Aut(A)$ be a $C_0(X)$-linear action of
$\RR^n$ on the separable continuous-trace $\cs$-algebra $A$ with spectrum $X$,
such that $\dim:X\to \ZZ^+$ is continuous. Let
$Y:=\specnp{(A\rtimes_{\alpha}\RR^n)}$ and let $p:Y\to X$
denote the canonical projection. Then
$\delta(A\rtimes_{\alpha}\RR^n)=p^*\delta(A)$, where
$\delta(A\rtimes_{\alpha}\RR^n)\in H^3(Y,\ZZ)$ and
$\delta(A)\in H^3(X,\ZZ)$ denote the
Dixmier-Douady invariants of $A$ and $A\rtimes_{\alpha}\RR^n$, respectively.
\end{thm}

For the proof we need
\begin{lem}\label{lem-Dix}
Let $p:Y\to X$ be a continuous map, $A$ a
continuous-trace \cs-algebra with spectrum $X$ and
$B$ a continuous-trace \cs-algebra with spectrum $Y$.
Suppose further that $B$ has the structure of a
$C_0(X)$-algebra via composition with $p$.
Then $A\otimes_XB$ is a continuous trace algebra
with Dixmier-Douady invariant
$\delta(A\otimes_XB)= p^*\delta(A)+\delta(B)$.
\end{lem}
\begin{proof} Recall that the pull-back $p^*A$ of $A$ along
$p$ is defined as $p^*A=A\otimes_XC_0(Y)$, where $C_0(Y)$ is given the
structure of a $C_0(X)$-algebra via the map
$C_0(X)\to C^b(Y);g\mapsto g\circ p$. On the other hand, we
may identify the $C_0(X)$-algebra $B$ with
$C_0(Y)\otimes_YB$, where the $C_0(X)$-action is given by the
$C_0(X)$-action on the first factor. Thus, by the associativity
of balanced tensor products, it follows that
$$p^*A\otimes_YB=(A\otimes_XC_0(Y))\otimes_YB=A\otimes_X(C_0(Y)\otimes_Y
B)=A\otimes_XB.$$
Since $\delta(p^*A)=p^*\delta(A)$ by \cite[Proposition 1.4]{rw},
the result follows from the multiplicativity of the Dixmier-Douady
invariant:
$\delta(p^*A\otimes_YB)=\delta(p^*A)+\delta(B)$ (cf., e.g.,
\cite[Prosition~2.2]{ckrw}).
\end{proof}

\begin{proof}[Proof of Theorem~\ref{Dix-Douady}]
First, since $\dim(x)\le n$ for all $x$, $X$ decomposes into
a finite disjoint
union
of open subsets such that $\dim$ is constant on each
subset.
Therefore, we may assume that $\dim$ is constant with value $k$, say, on all
of $X$.
Let $H$ be the representation group of $\RR^n$ and let
$$D_k:= \{[\om]\in H^2(\RR^n,\TT):\dim(\Sigma_{[\om]})=k\}.$$
Then it follows from \cite[Theorem 6.3.3]{ech4} that $D_k$
is a locally closed (and hence locally compact) subset of $H^2(\RR^n,\TT)$,
and that the restriction $\cs(H)_k:=\cs(H)_{D_k}$ of $\cs(H)$ to $D_k$ is a
continuous trace subquotient of $\cs(H)$. Thus, by Lipsman's and
Rosenberg's result \cite[Theorem 3.4]{lipros},
the Dixmier-Douady invariant of $\cs(H)_k$ is trivial.

Since we already assumed that $\dim$ has constant value $k$ on all of $X$,
it follows that $\varphi_{\alpha}$, and hence also the map
$f:X\to H^2(\RR^n,\TT)$ given by $ f(x)=\varphi_{\alpha}(x)^{-1}$
takes values in
$D_k$. Thus by Theorem~\ref{thm-rn} we get
$$A\rtimes_{\alpha}\RR^n\cong
A\otimes_Xf^*(\cs(H))=A\otimes_Xf^*(\cs(H)_k),$$
and $Y:=\specnp{(A\rtimes_{\alpha}\RR^n)}=f^*\widehat{H}_k$,
where $\widehat{H}_k$ denotes the restriction of the topological
bundle $\widehat{H}$ over $H^2(\RR^n,\TT)$ to $D_k$.
Since $Y=\specnp{f^*(\cs(H)_k)}$, it follows that
$f^*(\cs(H)_k)$ is a $C_0(Y)$-algebra in the canonical way and the
original $C_0(X)$-structure is induced via the projection $p:Y\to X$.
Thus by Lemma~\ref{lem-Dix} we have
$$\delta(A\rtimes_{\alpha}\RR^n)=\delta(A\otimes_Xf^*(\cs(H)_k))
=p^*\delta(A)+\delta(f^*(\cs(H)_k))=p^*\delta(A),$$
since $\delta(f^*(\cs(H)_k))=f^*\delta(\cs(H)_k)=0$.
\end{proof}

\begin{remark} Similar to the proof of \cite[Lemma 3.3]{lipros}
one can show that if $\alpha:\RR^n\to \Aut(A)$
is a $C_0(X)$-linear action of $\RR^n$ on the separable continuous trace
algebra $A$ with spectrum $X$, then any continuous trace subquotient
$B$ of $A\rtimes_{\alpha}\RR^n$ decomposes into a finite
direct sum of ideals such that all these ideals are subquotients of
some $A_{D_k}\rtimes_{\alpha}\RR^n$, where $D_k$ is a the locally closed
subset of
$X$ such that the dimension function is constantly equal to $k$ on $D_k$.
From this and the above result it follows that if $\delta(A)$ is trivial,
then any continuous trace subquotient of $A\rtimes_{\alpha}\RR^n$
has also trivial Dixmier-Douady invariant. We omit the details.
\end{remark}

%% file: cxcrossed.bbl
\def\mathcs{C^{\displaystyle *}}
  \def\cs{\ifmmode\mathcs\else$\mathcs$\fi}\def\noopsort#1{}
\ifx\undefined\bysame
\newcommand{\bysame}{\leavevmode\hbox to3em{\hrulefill}\,}
\fi
\begin{thebibliography}{10}

\bibitem{as}
Robert~J. Archbold and Donald W.~B. Somerset, {\em Quasi-standard
  {\cs}-algebras}, Math. Proc. Cambridge. Philos. Soc. {\bf 107} (1990),
  349--360.

\bibitem{bagpac}
Larry Baggett and Judith~A. Packer, {\em The primitive ideal space of two-step
  nilpotent group {\cs}-algebras}, J. Funct. Anal {\bf 124} (1994), 389--426.

\bibitem{blanchard2}
\'Etienne Blanchard, {\em Tensor products of {$C(X)$}-algebras over {$C(X)$}},
  Ast{\'e}risque {\bf 232} (1995), 81--92.

\bibitem{blanchard}
\bysame, {\em D\'eformations de {\cs}-alg\`ebres de {H}opf}, Bull. Soc. Math.
  France {\bf 124} (1996), 141--215.

\bibitem{cohen}
Paul~J. Cohen, {\em Factorization in group algebras}, Duke Math. J. {\bf 26}
  (1959), 199--205.

\bibitem{corgr}
Lawrence~J. Corwin and Frederick~P. Greenleaf, {\em Representations of
  nilpotent {L}ie groups and their applications, {P}art {I}: basic theory and
  examples}, Cambridge studies in advanced mathematics, vol.~18, Cambridge
  University Press, 1989.

\bibitem{ckrw}
David Crocker, Alexander Kumjian, Iain Raeburn, and Dana~P. Williams, {\em An
  equivariant {B}rauer group and actions of groups on {\cs}-algebras}, J.
  Funct. Anal. {\bf 146} (1997), 151--184.

\bibitem{ech4}
Siegfried Echterhoff, {\em Crossed products with continuous trace}, Mem. Amer.
  Math. Soc. {\bf 123} (1996), no.~586, 1--134.

\bibitem{echros}
Siegfried Echterhoff and Jonathan Rosenberg, {\em Fine structure of the
  {M}ackey machine for actions of abelian groups with constant {M}ackey
  obstuction}, Pacific J. Math. {\bf 170} (1995), 17--52.

\bibitem{ew2}
Siegfried Echterhoff and Dana~P. Williams, {\em Locally inner actions on
  {$C_0(X)$}-algebras}, preprint, June 1997.

\bibitem{fell-doran}
J.~M.~G. Fell and Robert Doran, {\em Representations of {$*$}-algebras, locally
  compact groups, and {B}anach {$*$}-algebraic bundles}, vol. {I\&II}, Academic
  Press, New York, 1988.

\bibitem{gior-mingo}
Thierry Giordano and James~A. Mingo, {\em Tensor products of {\cs}-algebras
  over abelian subalgebras}, J. London Math. Soc. {\bf 55} (1997), 170--180.

\bibitem{horr}
Steven Hurder, Dorte Olesen, Iain Raeburn, and Jonathan Rosenberg, {\em The
  {C}onnes spectrum for actions of abelian groups on continuous-trace
  algebras}, Ergod. Th. \& Dynam. Sys. {\bf 6} (1986), 541--560.

\bibitem{husemoller}
Dale Husemoller, {\em Fibre bundles}, Graduate Texts in Mathematics, vol.~20,
  Springer-Verlag, New York, 1975.

\bibitem{kw}
Eberhard Kirchberg and Simon Wassermann, {\em Operations on continuous bundles
  of {\cs}-algebras}, Math. Ann. {\bf 303} (1995), 677--697.

\bibitem{kmrw}
Alexander Kumjian, Paul~S. Muhly, Jean~N. Renault, and Dana~P. Williams, {\em
  The equivariant {B}rauer group of a locally compact groupoid}, preprint,
  1996.

\bibitem{lee2}
Ru~Ying Lee, {\em On the \cs-algebras of operator fields}, Indiana Univ. Math.
  J. {\bf 25} (1976), 303--314.

\bibitem{lipros}
Ronald Lipsman and Jonathan Rosenberg, {\em The behavior of fourier transforms
  for nilpotent {L}ie groups}, Trans. Amer. Math. Soc. {\bf 348} (1996),
  1031--1050.

\bibitem{moore3}
Calvin~C. Moore, {\em Group extensions and cohomology for locally compact
  groups. {III}}, Trans. Amer. Math. Soc. {\bf 221} ({\noopsort{d}}1976),
  1--33.

\bibitem{moore4}
\bysame, {\em Group extensions and cohomology for locally compact groups.
  {IV}}, Trans. Amer. Math. Soc. {\bf 221} ({\noopsort{e}}1976), 34--58.

\bibitem{may1}
May Nilsen, {\em {\cs}-bundles and {$C_0(X)$}-algebras}, Indiana Univ. Math. J.
  {\bf 45} (1996), 463--477.

\bibitem{judy-exp}
Judith~A. Packer, {\em Transformation group {\cs}-algebras: a selective
  survey}, {\cs}-algebras: 1943--1993: A fifty year Celebration (Robert~S.
  Doran, ed.), Contemp. Math., vol. 169, Amer. Math. Soc., 1994, pp.~182--217.

\bibitem{judymc}
\bysame, {\em {M}oore cohomology and central twisted crossed product
  {\cs}-algebras}, Canad. J. Math. {\bf 48} (1996), 159--174.

\bibitem{para1}
Judith~A. Packer and Iain Raeburn, {\em Twisted crossed products of
  {\cs}-algebras}, Math. Proc. Camb. Phil. Soc. {\bf 106} (1989), 293--311.

\bibitem{para3}
\bysame, {\em Twisted crossed products of {\cs}-algebras. {II}}, Math. Ann.
  {\bf 287} (1990), 595--612.

\bibitem{para2}
\bysame, {\em On the structure of twisted group {\cs}-algebras}, Trans. Amer.
  Math. Soc. {\bf 334} (1992), 685--718.

\bibitem{palais}
Richard~S. Palais, {\em On the existence of slices for actions of non-compact
  {L}ie groups}, Ann. Math. {\bf 73} (1961), 295--323.

\bibitem{pr2}
John Phillips and Iain Raeburn, {\em Crossed products by locally unitary
  automorphism groups and principal bundles}, J. Operator Theory {\bf 11}
  (1984), 215--241.

\bibitem{ra2}
Iain Raeburn, {\em Induced {\cs}-algebras and a symmetric imprimitivity
  theorem}, Math. Ann. {\bf 280} (1988), 369--387.

\bibitem{rr}
Iain Raeburn and Jonathan Rosenberg, {\em Crossed products of continuous-trace
  {\cs}-algebras by smooth actions}, Trans. Amer. Math. Soc. {\bf 305} (1988),
  1--45.

\bibitem{rw}
Iain Raeburn and Dana~P. Williams, {\em Pull-backs of {\cs}-algebras and
  crossed products by certain diagonal actions}, Trans. Amer. Math. Soc. {\bf
  287} (1985), 755--777.

\bibitem{rw2}
\bysame, {\em Crossed products by actions which are locally unitary on the
  stabilisers}, J. Funct. Anal. {\bf 81} (1988), 385--431.

\bibitem{90a}
\bysame, {\em Moore cohomology, principal bundles, and actions of groups on
  {\cs}-algebras}, Indiana U. Math. J. {\bf 40} (1991), 707--740.

\bibitem{90c}
\bysame, {\em Dixmier-{D}ouady classes of dynamical systems and crossed
  products}, Canad. J. Math. {\bf 45} (1993), 1032--1066.

\bibitem{rieffpr}
Marc~A. Rieffel, {\em Proper actions of groups on {\cs}-algebras}, Mappings of
  operator algebras (H.~Araki and R.~V. Kadison, eds.), Progr. Math., vol.~84,
  Birkhauser, Boston, 1988, Procceedings of the {J}apan-{U}.{S}. joint seminar,
  {U}niversity of Pennsylvania, pp.~141--182.

\bibitem{ros2}
Jonathan Rosenberg, {\em Some results on cohomology with {B}orel cochains, with
  applications to group actions on operator algebras}, Operator Theory:
  Advances and Applications {\bf 17} (1986), 301--330.

\bibitem{ros5}
\bysame, {\em {\cs}-algebras and {M}ackey's theory of group representations},
  {\cs}-algebras: 1943--1993: A fifty year Celebration (Robert~S. Doran, ed.),
  Contemp. Math., vol. 169, Amer. Math. Soc., 1994, pp.~150--181.

\bibitem{steenrod}
Norman~E. Steenrod, {\em Topology of fibre bundles}, Princeton University
  Press, New Jersey, 1951.

\bibitem{taktext}
Masamichi Takesaki, {\em Theory of operator algebras {I}}, Springer-Verlag, New
  York, 1979.

\bibitem{var}
V.~S. Varadarajan, {\em Geometry of quantum theory}, vol.~II, Van
  Nostrand-Reinhold, New York, 1970.

\bibitem{warner}
Frank~W. Warner, {\em Foundations of differentiable manifolds and {L}ie
  groups}, Scott, Foresman and Company, Glenview, Illinois, 1971.

\bibitem{wegge}
N.~E. Wegge-Olsen, {\em {K}-theory and {\cs}-algebras}, Oxford University
  Press, New York, 1993.

\end{thebibliography}
